\documentclass[11pt]{article}

\RequirePackage{fix-cm}
\usepackage{mathptmx}

\usepackage{latexsym}
\usepackage{marvosym}
\usepackage{amssymb}
\usepackage{amsfonts}
\usepackage{amsmath}
\usepackage{mathrsfs}
\usepackage{dsfont}    
\usepackage{bbold}     
\usepackage[english]{babel}
\usepackage{caption}
\usepackage{epsfig}
\usepackage{float}
\usepackage{subfigure}
\usepackage{psfrag}
\usepackage{graphicx}
\usepackage{epsfig}

\newcommand{\e}{\varepsilon}

\newcommand{\om}{\omega}
\newcommand{\g}{\gamma}

\newcommand{\ZZZ}{\mathds{Z}}

\newcommand{\NNN}{\mathds{N}}

\newcommand{\RRR}{\mathds{R}}
\newcommand{\TTT}{\mathds{T}}

\newcommand{\calK}{{\mathcal K}}
\newcommand{\st}{\scriptstyle}
\newcommand{\io}{\infty}

\def\v#1{\mathbf{#1}}
\newcommand{\pr}[2]{P(#1/#2)}
\newcommand{\ymax}{y_{\mathrm{max}}}
\newcommand{\Tmax}{T_{\mathrm{max}}}
\newcommand{\num}{\mathrm{num}}
\newcommand{\hem}{\mathrm{hem}}
\newcommand{\mNO}{\multicolumn{1}{|c|}{NO}}
\newcommand{\mNE}{\multicolumn{1}{|c|}{NE}}
\newcommand{\Del}{h}
\newcommand{\str}{\rule[-0.4ex]{0ex}{3ex}}

\begin{document}
\title{The high-order Euler method and the spin-orbit model\\
{\small A fast algorithm for solving differential equations with small,
smooth nonlinearity}}

\author
{\bf Michele Bartuccelli$^1$, Jonathan Deane$^1$, Guido Gentile$^2$
\vspace{2mm}
\\ \small 
$^1$ Department of Mathematics, University of Surrey, Guildford, GU2 7XH,
UK 
\\ \small
$^2$ Dipartimento di Matematica, Universit\`a di Roma Tre, Roma, I-00146,
Italy
\\ \small 
E-mail: m.bartuccelli@surrey.ac.uk, j.deane@surrey.ac.uk,
gentile@mat.uniroma3.it
}
\date{}
\maketitle

\begin{abstract}
We present an algorithm for the rapid numerical integration of smooth, time-periodic
differential equations with small nonlinearity, particularly suited to problems with
small dissipation. The emphasis is on speed without
compromising accuracy and we envisage applications in problems where
integration over long time scales is required; for instance, orbit
probability estimation via Monte Carlo simulation. We demonstrate the
effectiveness of our algorithm by applying it to the spin-orbit problem,
for which we have derived analytical results for comparison with those that we obtain
numerically. Among other tests, we carry out a careful comparison of our
numerical results with the analytically predicted set of periodic orbits 
that exists for given parameters. Further tests concern the long-term
behaviour of solutions moving towards the quasi-periodic attractor, and
capture probabilities for the periodic attractors computed 
from the formula of Goldreich and Peale.  We implement the algorithm in 
standard double precision arithmetic and show that this is adequate to obtain an
excellent measure of agreement between analytical predictions and the
proposed fast algorithm.
\end{abstract}

\section{Motivation}
\label{mot_sec}

In this paper, we discuss an algorithm for the rapid numerical solution of 
smooth, nonlinear, non-autonomous, time-periodic, dissipative 
differential equations, with reference to a particular example,
known as the spin-orbit equation. The spin-orbit ordinary differential
equation (ODE) describes the coupling, in the presence of tidal friction, between
the orbital and rotational motion of an ellipsoidal satellite orbiting a primary, and
many authors have studied it since the original work of~\cite{danby} and~\cite{GP};
see also~\cite{MD},~\cite{celbk} and~\cite{CL}.
In cases of interest, both the nonlinear and the dissipative terms are multiplied by
small parameters, and as the dissipation parameter decreases, the ODE possesses 
an ever-increasing number of co-existing periodic orbits, with the initial 
conditions selecting which one is observed. In this sense, the problem is 
not simple, despite the fact that the nonlinearity is small: small
dissipation coupled with small nonlinearity leads here to non-trivial dynamics.

One interesting application of the spin-orbit equation is as a model of
the orbit of Mercury, whose primary is considered to be the Sun; other
applications come to mind with the discovery of extra-solar planetary
systems. The orbit of Mercury appears to be unique in the solar
system, since it rotates three times on its own axis for every two
orbits of the Sun: all other regular satellites for
which we have data are in a one-to-one resonance with their primaries.
See for instance~\cite{noyelles} for a recent survey offering a new
perspective on the problem.

In order to estimate numerically the probability of capture of a satellite in a given
orbit, one possibility is to use a Monte Carlo approach, in which
the spin-orbit ODE is integrated forward in time, starting from many 
uniformly-distributed random initial conditions.  The time-asymptotic behaviour,
that is, the solution after any transient has decayed, is
determined for each of these initial conditions, and the probability of
capture by each of the possible attractors is thereby estimated. The challenges of 
this approach are (a) that realistic values of the dissipation parameter $\gamma$
are small, so transient times, which are $O(1/\gamma)$ --- see~\cite{BDG} for 
an argument in a similar case --- are long;
and (b), in order to obtain low-error estimates of capture probabilities, 
a large number $I$ of initial conditions must be considered: in fact,
the width of the 95\% confidence interval for the probabilities is proportional
to $I^{-1/2}$ --- see
equation~(\ref{conf}). In interesting cases, that is when $\gamma$ is small, 
(a) and (b) force one to carry out a large number of simulations of orbit 
dynamics, each one over a long time interval, which, using traditional 
numerical ODE solvers, requires prohibitively long computation times.
For a problem such as this, we therefore conclude that a fast ODE solving
algorithm is a necessity and not a luxury.

Many problems in mathematical physics boil down to solving an ODE for which
no closed-form solution exists. For simulations in such cases, there is no
alternative but to approximate solutions numerically. Also, the solutions
to nonlinear problems can
display sensitive dependence to initial conditions. This raises questions
as to how good a representation of what we casually refer to as
`the solution' to an initial value problem, is actually obtainable
numerically.  Contrast a finite precision numerical solution to the 
notional `true solution' --- one which is computed
to infinitely high precision, but the computation of which can be done
in a finite time. Clearly the latter is unattainable
with real computing hardware, with its finite memory and speed.
Hence, in practice, the best we can do is to use finite-precision,
usually double precision (typically 16--17 significant figures) algorithms,
to model, approximately, the true solution. 

Although software for arbitrary-precision arithmetic is available, we want
to show here what can be achieved using only double precision (with one
exception). The question then becomes: how might we construct a practical 
algorithm to approximate the true 
solution, using standard double precision arithmetic, while also bearing in mind
the need to obtain solutions quickly?

We describe in this paper an algorithm that speeds up the solution process by
a factor of at least 7 compared to `traditional' numerical ODE solvers, such
as standard algorithms like Runge-Kutta~\cite{Numrec,Asher}
and symplectic numerical methods, for instance the Yoshida 
algorithm~\cite[Appendix F]{yosh,celbk}.  The latter has been used to
solve the spin-orbit problem in the past, for example in~\cite{CC1}.
These algorithms and many more like them are general-purpose methods that
work for a wide variety of problems. By contrast, our algorithm is specific
to a particular problem, but, since it is set up by computer algebra, only
small changes need to be made to the set-up code in order to
modify it for a different problem; with this proviso, our algorithm is also
general-purpose. 

Our algorithm works well for problems like the spin-orbit ODE, for which we 
carry out a careful comparison of our numerical results with those
obtained analytically, via perturbation theory, as well as published
results on attractor probabilities, in order to validate our work. Setting
up the algorithm relies on computer algebra, and running the
algorithm at speed requires a low-level computer language; the interplay between these
two forms of computation is a theme in the paper.

The rest of the paper is organised as follows. The spin-orbit ODE is given
and a fast solution algorithm is described in Sect.~\ref{alg_sec}. Details
on setting up the algorithm and some practical data are given in
Sect.~\ref{setup_sec}, and verification is reported in
Sect.~\ref{ver_sec}. In Sect.~\ref{perf_sec}, we give details of the
speed and the robustness of the algorithm, and in Sect.~\ref{conc_sec}
we draw some conclusions.  The perturbation theory calculations which underpin the 
verifications are carried out in the Appendix, which also contains some
further supplementary material.

\section{The ODE and a fast solution algorithm}
\label{alg_sec}

We consider the spin-orbit ODE:
\begin{equation}
\begin{cases}
\dot  x = y , \\
\dot y = -\e \, G( x,t)-\gamma \alpha \left( y - \omega \right) .
\end{cases}
\label{ode}
\end{equation}
where  $\alpha, \e, \g, \omega >0$ and $x\in\TTT=\RRR/\pi\ZZZ$,
so that the phase space is $\TTT\times\RRR$.
Here $\e$ is a small parameter,
related to the asymmetry of the equatorial moments of inertia of the
satellite, and $e$ is the eccentricity of the orbit~\cite{GP,MD}. From here onwards
we set
$\dot x = y$. We follow~\cite{GP,CC1,CC2} in setting
$\omega = \nu(e) = \overline{N}(e)/\overline{L}(e)$; we also write
$\overline{L}(e) = \alpha$, where
$$\overline{L}(e) = \frac{1 + 3e^2 + 3e^4/8}{(1-e^2)^{9/2}}
\;\;\;\mbox{ and }\;\;\;
\overline{N}(e) = \frac{1 + 15e^2/2 + 45e^4/8 + 5e^6/16} {(1-e^2)^6}.$$
Furthermore, 
\begin{equation}
G(x,t)=\sum_{k\in\mathcal K} A_k(e)\sin(2x - kt),
\;\;\;\mbox{ where }\;\;\; \mathcal K = \{-3, -2, -1, 1, 2, 3, 4, 5, 6, 7\}
\label{gdef}
\end{equation}
and
\begin{eqnarray*}
A_{-3} = {\frac {81}{1280}}\,{e}^{5}
&\;\;\;\; &
A_{-2} = \frac{1}{24}\,{e}^{4}
\\
A_{-1} = {\frac {1}{48}}\,{e}^{3} + {\frac {11}{768}}\,{e}^{5}
&\;\;\;\; &
A_1 =  -\frac{1}{2}\,e + \frac{1}{16}\,{e}^{3} - {\frac {5}{384}}\,{e}^{5}
\\
A_2 = 1 - \frac{5}{2}\,{e}^{2} + {\frac {13}{16}}\,{e}^{4}
&\;\;\;\; &
A_3 = \frac{7}{2}\,e - \frac {123}{16}\,e^3 + {\frac {489}{128}}\,e^5
\\
A_4 =  {\frac {17}{2}}\,{e}^{2}-{\frac {115}{6}}\,{e}^{4}
&\;\;\;\; &
A_5 =  {\frac {845}{48}}\,{e}^{3}-{\frac {32525}{768}}\,{e}^{5}
\\
A_6 = {\frac {533}{16}}\,{e}^{4}
&\;\;\;\; &
A_7 = {\frac {228347}{3840}}\,{e}^{5}.
\\
\end{eqnarray*}
The expressions for $\overline{L}(e)$ and $\overline{N}(e)$ have been obtained
by averaging, and those for $A_k(e)$ have been derived by solving the Kepler
relations up to $O(e^6)$~\cite{CC1}, truncation at this order leading to
the neglect of all harmonics outside the set $\mathcal K$.

The dissipation model in equation~(\ref{ode}) is known as MacDonald's
tidal torque~\cite{macd,MD}. It has been widely studied
since the pioneering work of Goldreich and Peale~\cite{GP}, even though its
validity has recently been questioned; see for instance~\cite{noyelles}
and references therein, and also the comments at the end of Sect.~\ref{conc_sec}.
It should be noted that the probability of capture will be affected
by the choice of dissipation model.

The algorithm to solve~(\ref{ode}) that we propose in this paper is essentially
the usual Euler method, extended so that the series solution is computed to $O(\Del^N)$,
where $\Del$ is the timestep and $N\gg 1$. That is, we advance a solution by
one timestep via the truncated Taylor expansion
\begin{equation}
\v{x}(t_i) = \v{H}(\v{x}(t_{i-1}), t_{i-1}) = 
\v{x}(t_{i-1}) + \sum_{j=1}^N \frac{\Del^j}{j!}\, \v{f}_j(\v{x}(t_{i-1}), t_{i-1}),
\label{eul}
\end{equation}
where $t_i=t_0+i\Del$, ${\bf x}(t)=(x(t),y(t))$ and the functions ${\bf f}_{j}$
can be computed explicitly from the differential equation, which allows one to
compute, recursively, the derivatives of $x(t)$ and $y(t)$ of all orders at
$t=t_{i-1}$, in terms of the initial conditions, $x(t_{i-1})$, $y(t_{i-1})$.
The standard Euler method can be recovered by setting $N = 1$.

With a judicious choice of $N$ and $h$, we demonstrate 
that for our problem, one can use~(\ref{eul}) to compute solutions to the ODE in relatively
large, equal timesteps. Furthermore, the size of the timestep used is fixed throughout,
so the algorithm is not even adaptive. Such an approach might be thought to be
of limited practical use, but it is one purpose of this paper to show
that, for some problems, this is not the case. In particular, the
computational cost of solving an ODE using the proposed method turns out to be lower
than all other algorithms against which it was compared.

We draw a parallel here between this work and that of, for 
instance,~\cite{saari,chang}, in which a series approach is also used to solve 
ODEs. There are however important differences between the approach of
Chang and Corliss and ours: first, the series used by them are computed, 
numerically, at each 
timestep; and second, they use appropriate variations on the standard ratio test
for convergence, to estimate the size of each timestep --- so their method
is adaptive. By
contrast, in this work, the timesteps are fixed and all series required are
pre-computed and stored: this approach can significantly speed up the algorithm by
reducing the computational overheads.
Both methods are, however, essentially numerical analytical continuation.

In setting up our algorithm, we use computer algebra (CA) to generate code
in a low-level language (LL), once only for each set of parameters,
which computes the functions appearing on the right-hand 
side of equation~(\ref{eul}). This LL code is in turn compiled and executed 
in order to produce results. It might be thought that the
LL step can be omitted, and the CA program can be used to carry out the whole task. 
It can; this approach would lead to a significant decrease in speed however, since CA software 
is generally designed for algebraic manipulation and is not optimised for numerical computation.
As an example, consider the sum
\begin{equation}
S(n) = \sum_{i = 1}^n \frac{(i+1)(i+3)}{i(i+2)(i+4)(i+6)},
\;\;\;\mbox{ where }\;\;\;
\lim_{n\rightarrow\infty} S(n) = \frac{9}{32},
\label{sum}
\end{equation}
whose evaluation requires $6n-1$ addition and $5n$ multiplication/division operations, and
which we use later for timing purposes.
\footnote{In practice, we define 1 CPU-sec as the
time taken to evaluate $S(6\times 10^7)$: it happens to be the case that
the evaluation of $S(6\times 10^7)$ takes 1 second of CPU time on the
computer used to do most of the computations in this paper. Of course,
simply by timing the evaluation of $S(6\times 10^7)$ on another computer,
one can scale times given in this paper to correspond to times for that computer.}
The obvious experiment shows that numerical evaluation of $S(n)$, for $n$, say,
$10^6$, using 17 significant figures, takes about 260 times longer using CA compared with LL.
This increase in speed comes at a cost however: standard LL codes using
in-built mathematical operations, although relatively fast, will always carry
out arithmetic to fixed precision --- double precision is standard, which equates to
about 16--17 s.f. On the other hand, CA can in principle be used to evaluate numerical expressions
to any specified precision, the upper limit being set only by memory and time constraints.
This implied trade-off between accuracy and speed guides us in setting up
the algorithm in practice. The compromise we have to make is encapsulated in:
\begin{quote}
Find the smallest integer $N$ and the largest fixed timestep $\Del$, such that the pair of 
power series of degree $N$, which advance the solution $\v{x}(t)$ of~(\ref{ode}) 
from $t = i\Del$ to $t = (i+1)\Del$, using~(\ref{eul}), for all $i\in\NNN$, both do so 
to within a given tolerance.
\end{quote} 
Increasing $\Del$ increases speed, since larger timesteps are
used, and increasing $N$ and/or decreasing $\Del$ both increase accuracy
in principle, but the exact relationship between these parameters is not
straightforward, since rounding errors come into play. It is clear, though, that
since the differential equation~(\ref{ode}) is $2\pi$-periodic in $t$,
we need to find the smallest integer $M$, where $\Del =
2\pi/M$, such that a suitable error criterion is met for the finite set $i = 1, \ldots, M$, for
all initial conditions $\v{x}(0)$ in some subset $\mathcal Q$ of $\RRR^2$,
in order for it to be met for all $i\in\NNN$.

In order to quantify numerical error, we compare estimates of the state vector
$\v{x}(t) = (x(t), y(t))$ at a time $t = T_1$,
computed from the state vector at $t = T_0$, where $T_1 > T_0$,
with the computation being carried out in two ways: using a high-precision 
numerical ODE solver (which we denote with the subscript `num'), and our high-order
Euler method (which we label `hem').

Hence, the requirements of the computer algebra software are:
\begin{enumerate}
\item efficient series manipulation;
\item ability to translate arbitrary algebraic expressions into a
low-level language;
\item ability to carry out floating point arithmetic to any given
precision;
\label{ap}
\item a selection of algorithms for purely numerical solution of
differential equations.
\label{bs}
\end{enumerate}
Items~(\ref{ap}) and~(\ref{bs}) above are necessary for 
making error estimates. The numerical algorithm chosen was a Gear 
single-step extrapolation method using Bulirsch-Stoer rational 
extrapolation~\cite{Numrec}, which is
good for computing high-accuracy solutions to smooth problems. We make
the assumption that results produced by this algorithm, for $\Del\in[0,
2\pi]$, $t_0\in[0, 2\pi]$ and initial conditions in $\mathcal Q$,
are both accurate (that is, close to the true solution) and
precise (that is, correct to a large number of significant figures). In
fact, using 30 significant figures for computation, and relative and absolute 
error parameters of $10^{-20}$, we believe that numerical solutions
accurate to about 20 s.f.\ can be obtained, and it is against these that 
our algorithm is compared.
 
The computer algebra software Maple has all the necessary attributes and was used
for this work; the low-level language used was C.

The approach we adopt is partly experimental, in that we show that the power series 
we obtain meet the error criterion described in Sect.~\ref{setup_sec}, by comparing 
high-accuracy numerical 
solutions from CA with those produced by our algorithm, implemented in LL, and then using
the results to choose optimal values of $N$, the series truncation order, and $M$, the 
number of timesteps per period of $2\pi$.

In more detail, the computation of $\v{f}_i$ in~(\ref{eul}) is carried out as follows.
The method of Frobenius assumes that the solution to an ODE,
expanded about the point $t = t_0$, can be written as an infinite series,
so that $x(t) = \sum_{i = 0}^\infty a_i (t - t_0)^i$.
Substituting this into~(\ref{ode}) gives a recursion formula for $a_{i+1}$
in terms of $a_j,\, j = 0\ldots i$. Hence, given $a_0$ and $a_1$, which
correspond to the two initial conditions $x_0 = x(t_0)$
and $y_0 = \dot{x}(t_0)$, we can find $a_j$ for $j =2\ldots N$, where $N$
can in principle be as large as desired.

Since the ODE~(\ref{ode}) is nonlinear, so is the recursion formula, and the closed-form
expressions for $a_i(e, \varepsilon, \gamma, x_0, y_0, t_0)$, as polynomials in the
six arguments, quickly become large as $i$ increases. 
Hence, practical considerations, principally memory and computer time 
constraints, (a) force us to minimise the number of unevaluated parameters ---
we use the minimum, just two, $x_0$ and $y_0$, substituting numerical values for the 
others --- and (b) bound the value of $N$. For the specific case of the spin-orbit problem, it 
has been found to be feasible to use $N$ up to at least 20. This part of
the computation is carried out by CA.

Our ultimate goal is to estimate the relative areas of the basins of attraction of
each of the attractive periodic solutions to equation~(\ref{ode}), for given 
values of the parameters
$\varepsilon$ and $\gamma$. A Monte Carlo approach is one possible way to do
this. For the case at hand, this approach requires us first to compute
$\v{x}_j = (x_j, y_j) = \left(x(2j\pi), y(2j\pi)\right)$, $j = 1\ldots J$
for a sufficiently large $J$ that any transient behaviour has
effectively decayed away, and for a large number $I$ of uniformly-distributed random 
initial conditions $\v{x}_0 = (x_0, y_0)$ in a given set $\mathcal Q$. From now on, 
we drop the subscript 0 on
the initial conditions where this does not lead to confusion.  Clearly we need an 
efficient means for computing the Poincar\'{e} map
$\v{P}: \mathbb R^2 \mapsto \mathbb R^2$ generated by~(\ref{ode}), which is defined by
$\v{x}_{k+1} = \v{P}(\v{x}_k)$.  In practice, $\v{P}$ cannot be computed from the 
series solution in $M = 1$ step:
this would require $\Del = 2\pi/M = 2\pi$ in the series for $x(t)$ and $y(t)$, and this is 
certainly too large. Moveable singularities of the solution in complex-time
would prevent the series from converging for such a timestep. Hence, we split $\v{P}$
into $M$ `sub-maps' so that
$\v{P}(\v{x}) = \v{p}_M \circ \v{p}_{M-1} \circ \ldots \circ\v{p}_1 (\v{x})$,
where $\v{p}_i(\v{x}) = (X_i(\v{x}), Y_i(\v{x}))$, with $X_i$ advancing $x$
from $t = (i-1)\Del$ to $t = i\Del$ and $Y_i$ advancing $y$ over the same interval.
In terms of the function $\v{H}$ in equation~(\ref{eul}), we set $t_0 = 0$
and $\v{x} = \v{x}(t_{i-1})$, from which $\v{p}_i(\v{x}) =
\v{H}(\v{x}, t_{i-1})$.

With $\Del = 2\pi/M$, we have
\begin{equation}
X_i(\v{x}) = \sum_{j = 0}^N a_{i, j}(\v{x})\,\Del^j + O\left(\Del^{N+1}\right)
\;\;\; \mbox{ and }\;\;\;
Y_i(\v{x}) = \sum_{j = 0}^{N-1} (j+1)a_{i, j+1}(\v{x})\,\Del^j + O\left(\Del^{N}\right),
\label{XYdef}
\end{equation}
where $i = 1, \ldots, M$. Also, $\v{x} = (x(t_{i-1}), y(t_{i-1}))$ is the solution and 
its derivative at $t_{i-1} = (i-1)\Del$; and $a_{i, j}(\v{x})$ are polynomials
in $y,\, \cos 2x,\, \sin 2x$ if $j>0$, with an additional linear term in $x$ if $j = 0$.
We designate this algorithm the high-order Euler method (HEM).

In practice, the expressions for $X_i(\v{x})$ and $Y_i(\v{x})$, $i = 1, \ldots, M$,
are computed for particular numerical values of $e$, $\varepsilon, \gamma, N$ and 
$\Del$. The fact that the spin-orbit equation~(\ref{ode}) is also $\pi$-periodic 
in $x$ implies that the functions $X$ and $Y$, for fixed $M$ and $N$ and with numerical
values for $e$, $\varepsilon, \gamma$ and $\Del$, can be written in one of two forms. The
first of these is the Fourier form
\begin{eqnarray}
X_i(\v{x}) & = & x + A_{i, 0}(y) + \sum_{j = 1}^F \e^j\left[A_{i, j}(y)\cos 2jx 
+ B_{i, j}(y)\sin 2jx\right],\nonumber\\
Y_i(\v{x}) & = & C_{i, 0}(y) + \sum_{j = 1}^F \e^j\left[C_{i, j}(y)\cos 2jx +
D_{i, j}(y)\sin 2jx\right],
\label{fourform}
\end{eqnarray}
where $F$ is a positive integer and $A_{i, j}, \ldots D_{i, j}$ are
polynomials in $y$, $h$ and the parameters of the problem.
Both $F$ and the degree of the polynomials depend on our
accuracy requirements and on $N$; typically, we find $F\approx 3$ for $\e
\leq 10^{-3}$.
The fact that $A_{i, j}(y), \ldots D_{i, j}(y)$ always have a common
factor of $\e^j$ is explained in Appendix D.
The second (polynomial) form is equivalent to the Fourier form and is
\begin{equation}
X_i(\v{x}) = x + \sum_{i, j, k} \alpha_{i, j, k} \e^{j+k}y^i c^j s^k,\;\;\;\;
Y_i(\v{x}) = \sum_{i, j, k} \beta_{i, j, k}\e^{j+k} y^i c^j s^k,
\label{polyform}
\end{equation}
where $\alpha_{i, j, k},\, \beta_{i, j, k}$ are constants, and from here
onwards, we set $c = \cos 2x$, $s = \sin 2x$. In practice,
we use CA to compute $X_i$ and $Y_i,\, i = 1, \ldots, M$ in the polynomial
form, to convert these into Horner form~\cite{Numrec} for efficient evaluation, and 
then to translate the result
into LL. There turns out to be very little difference in the computational effort required
to evaluate these expressions in the Fourier and polynomial forms, and in
this work we choose the latter.

\section{Setting up the algorithm}
\label{setup_sec}

We now study a pair of cases in more detail.
Throughout this section, we let $\mathcal Q = [0, \pi]\times [0, \ymax]$ be
the set of initial conditions, with $\ymax = 5$. The first component of the initial
condition need only be in the range $0$ -- $\pi$  because the spin-orbit
equation is $\pi$-periodic in $x$. We also follow~\cite{CC1} in fixing
$e = 0.2056$, the value appropriate to Mercury, so that
$\omega\approx 1.25584$; $\e = 10^{-3}$; and $\g = 10^{-5}$ and $10^{-6}$,
giving $\gamma\alpha\approx 1.36937\times 10^{-5}$ and $1.36937\times 10^{-6}$ respectively.
All numerical computations in CA are carried out to 30 s.f.
We refer to these parameter values, with $\g$ excluded, as the default
parameters. The default value of $e$ and the other parameter values are chosen because we
can then compare our results using HEM directly with results published in~\cite{CC1},
which were obtained using a Yoshida symplectic integrator~\cite{yosh},~\cite[Appendix F]{celbk}.

\begin{table}[h]
\centering
\begin{tabular}{|c|c|c|c|c|c|} \hline
\str $N$ & $M$ & Total no.\ of & Total $+/\times$ ops.&
(Max $e_x$, Max $e_y$), $\times 10^{-14}$ & (Max $e_x$, Max $e_y$), $\times 10^{-14}$\\
& & terms in $\v{P}(\v{x})$ &(Horner form) & $\g = 10^{-5}$ & $\g = 10^{-6}$\\
\hline
18 & 18 & 6698 & 4597/4650 & 2582.6, 464.1  & 2586.9, 464.2 \\
   & 19 & 6914 & 4821/4866 & 734.2, 55.6  & 734.2, 55.7 \\
   & 20 & 7109 & 4994/5112 & 267.0, 21.1  & 267.0, 21.1 \\
   & 22 & 7433 & 5314/5400 & 46.0, 3.7  & 46.0,  3.6 \\
   & 25 & 7994 & 5872/5949 & 5.9, 0.54  & 5.4, 0.56\\
   & {\bf 28}	& {\bf 8396} & {\bf 6359/6455} & {\bf 4.1, 0.45}  & {\bf 4.4, 0.52}\\
   & 31 & 8780 & 6780/6807 & 3.8, 0.54  & 4.2, 0.45\\ \hline
19 & 22 & 7992 & 5685/5779 & 11.7, 0.98 & 10.0,  0.86\\
   & 25 & 8526 & 6270/6350 & 3.8, 0.47  & 4.1, 0.46\\
   & 28 & 8906 & 6750/6870 & 3.5, 0.51  & 4.4, 0.60\\ \hline
20 & 19 & 7967 & 5492/5608 & 44.6, 2.5  & 43.0, 2.8  \\
   & 20 & 8153 & 5715/5824 & 16.3, 1.0  & 14.0, 1.1  \\
   & 21 & 8346 & 5888/6067 & 6.9, 0.56  & 6.6, 0.63\\
   & 22 & 8512 & 6042/6170 & 4.6, 0.45  & 3.8, 0.50\\
   & 25 & 9042 & 6666/6741 & 4.3, 0.47  & 5.0, 0.50\\ \hline
\end{tabular}
\caption{How the number of arithmetical operations required to compute one iteration of
the Poincar\'{e} map, and the approximate maximum error obtained
when using the high-order Euler Method, vary with $N$ and $M$. The maximum
error is an estimate of $\max_{\v{x}_0\in\mathcal Q} (e_x(\v{x}_0),
e_y(\v{x}_0))$.}
\label{table1}
\end{table}

We first use CA to set up the functions $\v{p}_i(\v{x})$ and then translate
them into LL.  {\em A priori}, we have no 
idea what values of $M$ and $N$ to choose, and a compromise between high
accuracy, which tends to increase $M$ and $N$, and speed of the HEM,
which increases with decreasing $M$, must be found. Additionally, finite computer
memory puts a bound on $N$, since the expressions for $a_{i, j}(\v{x})$
in equation~(\ref{XYdef})
grow rapidly in size with $j$. Furthermore, the fact that these expressions will
eventually be evaluated using finite-precision arithmetic means that
increasing $N$ and $M$ too much can result in a {\em less} accurate
approximation to the Poincar\'{e} map, owing to the fact that more
operations are required to evaluate the expressions, potentially leading
to increased rounding errors.

We define our measure of error as follows. Letting $\v{x}_0 = \v{x}(t_0)$, we define the error
vector $\v{e}(\v{x}_0) = (e_x, e_y)$ by
\begin{equation}
e_x = \left|x_\num(\v{x}_0, t_0+2\pi) - x_\hem(\v{x}_0, t_0+2\pi)\right|,\;\;
e_y = \left|y_\num(\v{x}_0, t_0+2\pi) - y_\hem(\v{x}_0, t_0+2\pi)\right|.
\label{err}
\end{equation}

In practice, we estimate the maximum values of $e_x(\v{x}_0)$ and
$e_y(\v{x}_0)$, with $t_0 = 0$, as $\v{x}_0$ ranges over a grid of uniformly-spaced
points in $\mathcal Q$. The points used are 
$\{\v{x}_0 = (i\,\Delta x, j\,\Delta y),\, i, j = 0\ldots L\}$ with
$\Delta x= \pi/L$ and $\Delta y = \ymax/L$, and $L = 25$.

We now establish good values of $M$ and $N$. Table~\ref{table1} gives data to
guide the choice of values that represents a compromise between accuracy and speed.
The total number of terms and operation count data, which are almost the
same for both $\g = 10^{-5}$ and $10^{-6}$, are given to enable us to judge the
relative speed, and the $e_x, e_y$ values indicate the accuracy.

The main point to note is that for fixed $N$, the maximum error varies
little with $M$ for $M \geq M_{\mathrm{crit}}$, but for $M <M_{\mathrm{crit}}$,
the error increases rapidly: there is a `knee' in the error curve at $M =
M_{\mathrm{crit}}$. From several possible candidates, we choose
$N = 18$ and $M = 28$, which represents a good speed/accuracy compromise both for 
$\g = 10^{-5}$ and $10^{-6}$.

\begin{table}[h]
\centering
\begin{tabular}{|l|l|l|} \hline
\multicolumn{1}{|c} {\raisebox{2.2ex}{ }Description} &
\multicolumn{1}{|c} {$\g = 10^{-5}$} &
\multicolumn{1}{|c|} {$\g = 10^{-6}$}\\ \hline
1. Total CPU time for computing $X_i(\v{x})$, $Y_i(\v{x})$,&&\\
$i = 1, \ldots M$, using CA. Without error check. & 255 CPU-sec & 262 CPU-sec\\ \hline
2. As 1., but with error check. & 3143 CPU-sec & 2550 CPU-sec\\ \hline
3. Total no.\ of terms in $\v{P}(\v{x})$, before/after pruning &\multicolumn{2}{c|} {20578/8396} \\ \hline
4. Number of $+/\times$ operations (pruned, not Horner) & \multicolumn{2}{c|} {6349/55550}\\ \hline
5. Number of $+/\times$ operations (pruned, Horner form) & \multicolumn{2}{c|} {6359/6455}\\ \hline
6. $\mbox{[Maximum value of $(e_x, e_y)$]}\times 10^{-14}$ & $(4.1, 0.45)$ & $(4.4, 0.52)$\\ \hline
\end{tabular}
\caption{Data on the computer algebra set-up of the low-level language code to compute 
the Poincar\'{e} map.}
\label{table2}
\end{table}

Table~\ref{table2} gives some data on setting up $\v{P}(\v{x})$ with $N = 18$ and $M = 28$.
Points to note, with numbers in the list corresponding to line numbers in
the table, are:

\begin{enumerate}
\item The timings are given here in units of seconds of CPU time on a
particular computer. For comparison, the LL-evaluation of $S(6\times 10^7)$
defined in equation~(\ref{sum}), using the same computer, takes about 1 CPU-sec.
\item The error check is an estimate of the maximum values of $e_x(\v{x}_0)$ and
$e_y(\v{x}_0)$, defined in equation~(\ref{err}), with $t_0 = 0$,
as $\v{x}_0$ ranges over a grid of uniformly-spaced points in $\mathcal Q$,
as defined immediately following equation~(\ref{err}). We set $L = 25$, so
that $\Delta x= \pi/25$ and $\Delta y = 0.2$.
\item The total number of terms in the expressions for $X_i(\v{x})$, $Y_i(\v{x})$, $i
= 1,\ldots, M$ in the polynomial form --- equation~(\ref{polyform}) --- is
given in row~3 in the table. A `term' is a single product of the form
$\alpha_{i,j,k} y^i c^j s^k$ appearing in equation~(\ref{polyform}).

`Pruning' is a way of cutting down the number of terms retained by removing the
negligible ones. Specifically, any terms for which $|\alpha_{i,j,k}\,\ymax^i| < \Tmax$
\footnote{In fact this is an overestimate of the maximum value of a
given term: taking into account the powers of $s = \sin x$ and $c = \cos x$, the maximum
value of the term should be multiplied by $[j/(j+k)]^{j/2} [k/(j+k)]^{k/2}$,
which is $\max_{x\in\mathbb R}\cos^j x\sin^k x$ and is of order 1 for
relevant values of $j, k$.}
are deleted, with $\Tmax = 10^{-18}$.
This value of $\Tmax$ was chosen because the final expressions will be
computed in in LL using double precision arithmetic (equivalent to about 17 s.f.).
Pruning with $\Tmax = 10^{-18}$
reduces the number of terms in $\v{P}(\v{x})$ by about a half.
\item This row gives a measure of the computational cost --- the total number of
addition and multiplication operations --- of evaluating $\v{P}(\v{x})$ in
the polynomial form. Evaluating a term $y^k$ is assumed to take $k-1$
multiplications.
\item The total number of multiplication operations is reduced by a factor
of about eight when the expressions are converted to Horner form.
\item For the given parameters, the maximum difference between one iteration of
the Poincar\'{e} map computed using (a) HEM and (b) a high-precision numerical ODE
solver, is of order $10^{-14}$. The values given are an approximation to
$\max_{\v{x}_0\in\mathcal Q} \left(e_x(\v{x}_0), e_y(\v{x}_0)\right)$, as defined in 
equation~(\ref{err}) with $t_0 = 0$.
\end{enumerate}

\section{Verification}
\label{ver_sec}

We now verify the HEM in three ways. In the first of these, we check that the 
set of periodic and quasi-periodic solutions obtained numerically via the HEM
corresponds with those that can be proved to exist analytically, for
example by perturbation theory. In the
second verification, the probability of obtaining the different attractors
is estimated, again using the HEM, and these probabilities are compared
with the results published in~\cite{CC1} --- here, we are of course
comparing one numerical algorithm (HEM), against another (a Yoshida
symplectic integrator). We also compare attractor probabilities given by the formula
of Goldreich and Peale~\cite{GP} with those obtained from HEM.
In the third verification, we compute
$\omega'$ (defined in the Appendix and discussed below) numerically, using the HEM, 
and compare this value with that given by perturbation theory.

All the verifications relate to the default parameters, which are listed at the start of
Sect.~\ref{setup_sec}. The probability computations using HEM are carried out in the
standard way: $I$ uniformly-distributed random initial conditions in $\mathcal Q$ are
selected and the Poincar\'{e} map is iterated $n_{\mathrm{pre}}$ times, starting
from each one. Since the transient time is $O(1/\g)$~\cite{BDG},
we use $n_{\mathrm{pre}} = 10^6$ for $\g = 10^{-5}$ and
$n_{\mathrm{pre}} = 5\times 10^7$ for $\g = 10^{-6}$. For the other values of
$\g$ used later in the paper, we also choose $n_{\mathrm{pre}} = m/\g$, with
$m\approx 10$ being chosen such that, after integrating for a time $2\pi
n_{\mathrm{pre}}$, any transients have decayed to the point where the solution
can be identified.

\subsection{Which attractors exist?}

As shown in the Appendix, for $\g = 10^{-5}$, the quasi-periodic solution
and periodic solutions with
$p/q = 1/2,\, 1/1,\, 5/4,\, 3/2,\, 2/1,\, 5/2$ and $3/1$
exist according to a second order analysis, and no others.
For $\g = 10^{-6}$, these solutions remain, and additional solutions with $p/q
= 3/4,\, 7/4$ and $7/2$ also exist. All of these, and only these solutions are observed 
when using the HEM,
although only two out of the $I = 32000$ random initial conditions were
attracted to the $p/q = 3/4$ solution. Since the threshold for this solution is
$\g = 1.058\times 10^{-6}$, we would expect the probability of observing it 
to be very small. \textit{A posteriori} we expect that higher order
periodic solutions do not arise for the chosen values of the parameters ---
or, at worst, are irrelevant.

\subsection{The quasi-periodic solution}

A quasi-periodic solution to~(\ref{ode}) can also exist, as discussed
in the Appendix. This solution has a mean growth rate, $\omega'$, which, for
small $\varepsilon$, is close to $\omega$.  Equation~(\ref{eq:3.1}) implies 
a formula for estimating $\omega'$, which is
\begin{equation}
\omega' = \lim_{t\rightarrow\infty}\left[\frac{x(t) - x(0)}{t}\right]
 = \lim_{n\rightarrow\infty}\left[\frac{x(2\pi n) - x(0)}{2\pi n}\right].
\label{omp_comp}
\end{equation}
The second version is appropriate here, since we use the HEM to compute
iterations of the Poincar\'{e} map, and so only have access to values of
$(x(t), \dot{x}(t))$ at $t = 2n\pi$, $n = 0, 1, 2,\ldots$. 

\begin{figure}[!ht]
\centering 
\includegraphics*[width=3in]{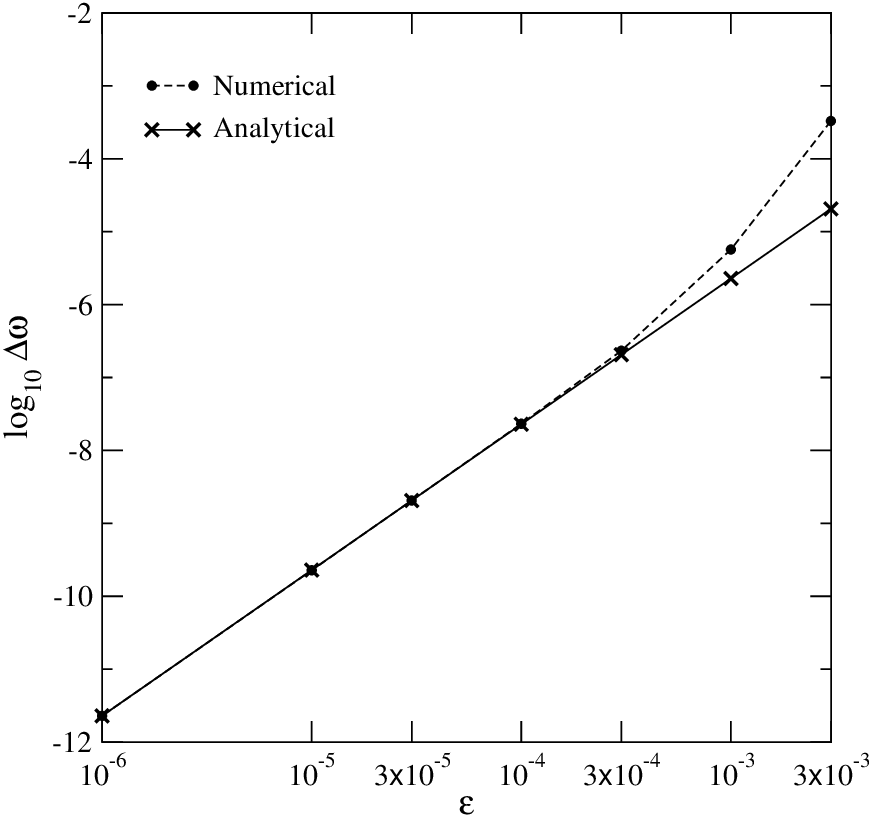}
\caption{Comparison of analytical and numerical computations of
$\Delta\omega = \omega - \omega'$ for various values of $\varepsilon$. For
small $\varepsilon$, the values of $\Delta\omega$ are seen to agree,
thereby validating the HEM, which was used to produce the numerical
results.}
\label{qp_comp}
\end{figure} 

Starting from equation~(\ref{omp}) in the Appendix, we have that 
$$\omega' = \omega - \varepsilon^2\mu^{(2)}(\omega') + O(\varepsilon^3)
= \omega - \varepsilon^2\mu^{(2)}(\omega) + O(\varepsilon^3),$$
where we have used the fact that,
since $\om'$ and $\om$ differ by $O(\varepsilon^2)$, replacing
$\mu^{(2)}(\om')$ with $\mu^{(2)}(\om)$ makes a difference
$O(\varepsilon^4)$, which can be neglected. From equation~(\ref{eq:3.14})
we compute $\mu^{(2)}(\om) = 2.284502$. In Fig.~\ref{qp_comp} we
plot the analytical estimate of $\Delta\om = \om - \om' =
\varepsilon^2\mu^{(2)}(\omega)$ and the numerical estimate
from~(\ref{omp_comp}), using LL, for $n = 10^8$, $\gamma = 10^{-5}$ and
$\varepsilon\in\{10^{-5}, 3\times 10^{-5}, 10^{-4}, 3\times 10^{-4}, 10^{-3}, 3\times
10^{-3}\}$. Additionally, we estimate $\Delta\omega$ for $\varepsilon = 10^{-6}$,
but in this case, double precision arithmetic is inadequate --- this the single
case, referred to in Sect.~1, where we do not use double precision.
Details of this computation are given in Sect.~C of the Appendix.

This is very different kind of test of the HEM compared to that
described in the previous section. Here, we check that the long-term
average rate of increase of $x$ implied by equation~(\ref{omp_comp}) is as
predicted by the analytical computation.

\subsection{Estimated attractor probabilities}

We estimate the probability $\pr{p}{q}$ that integrating forward in time from a 
randomly-selected initial condition $\v{x}\in\mathcal Q$ leads to a
period-$p/q$ orbit.  If several periodic orbits with a given $p$, $q$
exist, then their combined probability is computed.

\begin{table}[!ht]
\centering
\begin{tabular}{|l|l|l|l|l|} \hline
\multicolumn{1}{|c} {     } & \multicolumn{4}{|c|} {Probability, $\pr{p}{q}$, \%} \\ \cline{2-5}
$p/q$	& \multicolumn{2}{|c} {\str $\g = 10^{-5}$} & \multicolumn{2}{|c|} {$\g = 10^{-6}$} \\ \cline{2-5}
	& \multicolumn{1}{c|} {\str From CC} & \multicolumn{1}{c|}{This work} &
          \multicolumn{1}{c|}{From CC} & \multicolumn{1}{c|}{This work} \\ \hline
1/2       & \mNO          & $0.50\pm 0.08$   & \mNO          & $0.62\pm 0.09$\\
3/4       & \mNE          & \mNE              & \mNO          & $0.0063(\pm 0.009)$ \\
1/1       & $4.7\pm 1.3$ & $4.58\pm 0.23$   & $4.6\pm 1.3$  & $4.77\pm 0.23$\\
5/4       & $8.4\pm 1.7$ & $7.31\pm 0.29$   & $5.1\pm 1.4$  & $7.50\pm 0.29$\\
$\omega'$ & $69.8\pm 2.9$ & $71.65\pm 0.49$ & $73.4\pm 2.7$ & $70.22\pm 0.50$\\
3/2       & $12.6\pm 2.1$ & $12.05\pm 0.36$ & $14.0\pm 2.2$ & $11.94\pm 0.36$ \\
7/4       & \mNE           & \mNE             & \mNO           & $ 0.094\pm 0.03$\\
2/1       & $2.9\pm 1.0$ & $2.72\pm 0.18$   & $2.5\pm 0.97$  & $3.01\pm 0.19$\\
5/2       & $1.1\pm 0.7$ & $0.97\pm 0.11$   & $0.2(\pm 0.28)$ & $1.13\pm 0.12$\\
3/1       & $0.5\pm 0.4$ & $0.22\pm 0.05$   & $0.2(\pm 0.28)$ & $0.48\pm 0.08$\\
7/2       & \mNE          & \mNE              & \mNO            & $0.24\pm 0.05$\\ \hline
\end{tabular}
\caption{Attractor probabilities with their 95\% confidence intervals,
 determined using $I = 1000$ points, taken from CC \cite{CC1};
and 32000 points (this work). NO: attractor exists but was not observed; NE: 
attractor non-existent for these
parameters. The 95\% confidence interval in parentheses is not reliable since
for this case, $\widehat{p}I < 5$. The Poincar\'{e} map was
iterated a total of about $1.6\times 10^{12}$ times to produce the probability data
for $\g = 10^{-6}$.}
\label{vertab}
\end{table}

We also compute the 95\% confidence interval for these probabilities, using
the formula for the standard error of a proportion~\cite{Wal}.
This states that if a number $I$ of initial conditions is considered;
$\widehat{p}$ is the number of those initial conditions that end up on a given
attractor $A$, divided by $I$; and $Z_{c/2}$ is defined by 
$$\frac{1}{\sqrt{2\pi}}\int_{-Z_{c/2}}^{Z_{c/2}} e^{-z^2/2}\, dz = c, 
\;\;\; \mbox{ where $c\in[0, 1]$}; $$
then a $c\times 100\%$ confidence interval for the actual proportion $p$ of 
initial conditions going to $A$ is
\begin{equation}
p\in \left[\widehat{p}-Z_{c/2}\sqrt{\frac{\widehat{p}(1-\widehat{p})}{I}},\;
\widehat{p}+Z_{c/2}\sqrt{\frac{\widehat{p}(1-\widehat{p})}{I}}\right].
\label{conf}
\end{equation}
This estimate is reliable provided that $I\widehat{p}\geq 5$~\cite{Wal}.
Setting $c = 0.95$ corresponds to a 95\% confidence interval and gives $Z_{0.475}\approx 1.96$.
Clearly, the width of the confidence interval is proportional to $I^{-\frac{1}{2}}$,
as stated in Sect.~1. 

Note that the simulations reported in~\cite{CC1} do not find all possible 
periodic orbits.
Take the case $p/q = 1/2$ for $\g = 10^{-5}$. From Table~\ref{vertab}, we
have $P(1/2)\in [4.2\times 10^{-3},\, 5.8\times 10^{-3}]$ with 95\% confidence. 
From the binomial distribution, one can compute that the probability of this
orbit not being observed at all in 1000 simulations is less than $0.015$.

The periodic orbit probabilities for
$\g = 10^{-5}$ and $10^{-6}$ are given in Table~\ref{vertab}.
Note that the results for $\g = 10^{-6}$ in~\cite{CC1} were obtained by polynomial 
extrapolation from larger $\g$ values and the 95\% confidence intervals, added by 
us, were computed assuming that $I = 1000$.
Extrapolation in a case like this can be risky, because when taking smaller
values of $\gamma$ the orbit probabilities do not increase indefinitely, but tend
to settle around a constant value; this has been observed numerically in~\cite{BDG}
for a system with cubic nonlinearity, but we believe this to be a general phenomenon.
In our case, $\e = 10^{-3}$ with $e=0.2056$, the value of $\gamma$ where this 
appears to happen is around $\gamma = 10^{-5}$, and the constant value of
$\pr{3}{2}$ for $\gamma < 10^{-5}$ is about 12\%.

There is a further check that we can carry out, based on the formula of
Goldreich and Peale~\cite{GP}. Using an averaging technique, this formula
approximates the probability of capture in a particular $p:q$ resonance,
with $q = 2$, as follows:
$$P_{GP}(p) = \frac{2}{1 + \frac{\pi(p/2 - \omega)}{2\sqrt{2\e A_p(e)}}},$$
where $A_p(e)$ is defined straight after equation~(\ref{gdef}).
The formula can be seen to be $\gamma$-independent, but, for small enough
$\gamma$, gives probability estimates in good agreement with those given by
HEM, as shown in Table~\ref{GPtab}.

\begin{table}[!ht]
\centering
\begin{tabular}{|l|l|c|c|c|c|c|} \hline
& & \multicolumn{5}{|c|} {\str Probability, $\pr{3}{2}$, \%}\\ \cline{3-7}
\str $\epsilon$ & $\gamma$ &  & 
\multicolumn{2}{|c|}{\str $y\in[1.5, 2]$} & \multicolumn{2}{|c|}{$y\in[1.5, 5]$}\\ \cline{4-7}
& & \raisebox{1.2ex}{G \& P} &
HEM, all terms & HEM, $A_3$ only & HEM, all terms & HEM, $A_3$ only\\ \hline
\str                 & $10^{-8}$ & 7.70 & $9.17 \pm 0.48$ &$9.56 \pm 0.44$ & $7.82\pm 0.41$
&$8.34 \pm 0.42$\\
                     & $10^{-7}$ & 7.70 & $9.84 \pm 0.48$ & $9.16\pm 0.33$ & $7.43\pm 0.47$ & $7.84\pm 0.38$\\
$1.8\times 10^{-4}$  & $10^{-6}$ & 7.70 & $9.28 \pm 0.37$ & $9.49\pm 0.45$ & $7.92\pm 0.34$ & $7.88\pm 0.42$\\
                     & $10^{-5}$ & 7.70 & $9.02 \pm 0.36$ & $8.77\pm 0.36$ & $7.48\pm 0.37$ & $7.48\pm 0.33$\\
                     & $10^{-4}$ & 7.70 & $5.73\pm 0.29$  & $5.61\pm 0.29$ & $4.61\pm 0.27$ &$4.62\pm 0.27$\\ \hline
\str                 & $10^{-7}$ & 17.24 & $20.3 \pm 0.47$ & & $16.0\pm 0.43$ &\\
                     & $10^{-6}$ & 17.24 & $20.6 \pm 0.42$ & & $16.6\pm 0.27$ &\\
$1.0\times 10^{-3}$  & $10^{-5}$ & 17.24 & $20.2 \pm 0.28$ & & $16.3\pm 0.47$ &\\
                     & $10^{-4}$ & 17.24 & $19.0 \pm 0.34$ & & $15.8\pm 0.46$ & \\
                     & $10^{-3}$ & 17.24 & $8.59 \pm 0.36$ & & $6.80\pm 0.32$ & \\ \hline
\end{tabular}
\caption{Comparison of the probability of capture by the 3:2 orbit,
for $x\in[0, \pi], y\in [1.5, 2]$ and $y\in [1.5, 5]$, as computed by G\&P, 
the formula of Goldreich and Peale~\cite{GP} and also by the high-order Euler 
method, HEM. The $\pm$
quantities after the HEM probabilities are the width of the 95\% confidence
interval. These are given to two significant figures so that, $I$, the number of
initial conditions used, can be deduced if required, by using equation~(\ref{conf}).}
\label{GPtab}
\end{table}

The Goldreich and Peale formula is obtained under a series of
approximations, one of which consists in assuming that the solution is close
to a given resonance.  That is, the formula computes the probability of 
capture for a solution passing near the given resonance and the possibility
that the solution is captured by other resonances is neglected. Since, on physical
grounds, we are interested in trajectories in which the speed of rotation
decreases with time, it might be expected that the best choice would be to
take the initial data above the resonance 3:2 and below the next higher
resonance, {\em i.e.} $y\in [1.5, 2]$. However, as Table~\ref{GPtab} shows, the 
Goldreich and Peale formula better describes the behaviour of trajectories
starting in the full phase space above the resonance. 
\footnote{We explicitly consider data with $y\leq 5$, as in previous simulations,
and we have checked numerically that the probabilities do
not change appreciably when the initial velocity is further increased.}
The presence of the other resonances apparently does not affect
the probability of capture by the 3:2 resonance. We also considered a
modified model of the form~(\ref{ode}), where only the harmonic with $k=3$ is kept
in $G(x, t)$. Here we found that the probability of capture in the 3:2
resonance is essentially the same as for the full spin-orbit model;
apparently the basins of attraction of the other resonances are formed at
the expense of the basin of attraction of the quasi-periodic attractor,
leaving that of the 3:2 resonance unaffected.

\section{Performance of the high-order Euler method}
\label{perf_sec}

\subsection{Speed}
\label{speedsec}

The HEM was developed as a fast numerical solver 
for the spin-orbit ODE and problems like it. Hence, we now compare the timings 
for solving~(\ref{ode}) using HEM, with those from two other numerical methods. 
We choose an explicit Runge-Kutta method due to Dormand and Price, as described
in~\cite{hnw}, and an adaptive Taylor series method due to Jorba and
Zou (TSM)~\cite{jorba}.\footnote{Codes to implement these methods are available
for download. The Runge-Kutta code used here, DOP853, is available at {\tt
http://www.unige.ch/\textasciitilde hairer/software.html} and the Taylor series code can
be found at {\tt http://www.maia.ub.edu/\textasciitilde angel/taylor/}.}

\begin{table}[!ht]
\centering
\begin{tabular}{|l|l|l|l|l|} \hline
\multicolumn{1}{|c}{Parameters} &\multicolumn{1}{|c} {Tolerance} & 
\multicolumn{1}{|c|} {HEM time,} &\multicolumn{1}{|c|} {Ratio (DOP853)} &
\multicolumn{1}{|c|} {Ratio (Taylor} \\ 
\multicolumn{1}{|c}{ } & \multicolumn{1}{|c} { } & 
\multicolumn{1}{|c|} {CPU-sec} &\multicolumn{1}{|c|} { } &
\multicolumn{1}{|c|} {series method)} \\ \hline
$\str \varepsilon = 1.2\times 10^{-4}$, $\g = 10^{-7}$     & $3.6\times 10^{-15}$ & 17.22 & 21.1 & 11.7\\ \hline
$\str \varepsilon =1.8\times 10^{-4}$, $\g =10^{-8}$       & $4.4\times 10^{-15}$ & 18.86 & 19.0 & 11.0\\ \hline
$\str \varepsilon = 10^{-3}$, $\g = 5\times 10^{-6}$ & $4.5\times 10^{-14}$ & 19.04 & 14.2 & 7.97\\ \hline
$\str \varepsilon = 10^{-3}$, $\g = 10^{-6}$       & $2.1\times 10^{-14}$ & 20.93 & 16.3 & 9.05\\ \hline
$\str \varepsilon = 3\times 10^{-3}$, $\g = 10^{-5}$       & $4.2\times 10^{-14}$ & 24.44 & 12.6 & 6.92\\ \hline
\end{tabular}
\caption{Comparison of timings for the high-order Euler method (HEM) versus
a Runge-Kutta code (DOP853) and a Taylor series method. The timings are
the mean from three computations, in each of which the Poincar\'{e} map was 
iterated 50 000 times starting from each of 50 random initial conditions --- 
hence $2.5\times 10^6$ iterations for each computation. For HEM, the actual
time in CPU-sec is given; the last two columns give the {\it ratio} of the
time taken by the named algorithm to the time taken by HEM.}
\label{speedtab}
\end{table}

The results are summarised in Table~\ref{speedtab}, in which we compare the
data on the time taken for each of the algorithms to perform $2.5\times 10^6$
iterations of the Poincar\'{e} map. Specifically, if we define the test
problem as `iterate the Poincar\'{e} map 50 000 times starting from each of
50 random initial conditions in $\mathcal Q$', then the time used to
produce Table~\ref{speedtab} is the mean of the time taken to run the test
problem three times, using a different set of initial conditions on each occasion.
Note that the figure in the `HEM time' column is a number of CPU-sec,
where 1 CPU-sec is the time taken to compute $S(6\times 10^7)$,
defined in equation~(\ref{sum}). By contrast, the figures in the two `Ratio'
columns are the ratios of the times taken by the named algorithms to the time
taken by HEM.

We have been at pains to make the comparisons as fair as possible, which is
why the tolerance is different for each of the five sets of parameters: the
values chosen correspond closely to the estimated tolerance in the HEM
method for those parameters. We take this precaution because the time taken
by both Runge-Kutta and TSM depends sensitively on the value of tolerance used.

It can be seen from Table~\ref{speedtab} that HEM outpaces both the algorithms
against which it has been tested by a factor of at least 6.9:1, with the 
factor depending on the parameters. It is noteworthy that the time taken by
HEM is not strongly correlated with the value of $\gamma$, with, in
particular, the smaller (and physically more interesting)
values of dissipation not significantly slowing down the computation: in
fact, there is evidence that smaller values of $\gamma$ lead to a relative
\textit{increase} in speed of HEM compared to the other algorithms.

It may be thought surprising that HEM 
is noticeably faster than TSM, since both are based on analytical 
continuation. A quick experiment for the first set of parameters in
Table~\ref{speedtab} shows that, typically, the TSM takes about 15.5 timesteps
to advance a solution of the spin-orbit ODE by a time $2\pi$, which is
comparable with $M$, the number used in HEM ($M\approx 20-30$). 
Hence, the likely reason for the difference in speed is that,
since TSM is an adaptive algorithm, the Taylor series for the solution must 
be re-computed at every timestep.
This results in a larger computational overhead compared to HEM, where the
series are computed once only, saved, and then merely evaluated as required
in order to compute the Poincar\'{e} map.

\subsection{Robustness}

We now give some results that illustrate the robustness of the computation of
attractor probabilities using HEM. We deliberately choose sub-optimal
values of $M$ that result in higher maximum values of $(e_x, e_y)$, the
absolute error per iteration of the Poincar\'{e} map. The
data are given in Table~\ref{fuzztab}, from which it can be concluded that
the probabilities computed with all three values of $M$ agree at
the 95\% confidence level in the cases considered.

\begin{table}[!ht]
\centering
\begin{tabular}{|l|l|l|l|l|} \hline
\multicolumn{1}{|c} {     } & \multicolumn{4}{|c|} {Probability, $\pr{p}{q}$, \%} \\ \cline{2-5}
$p/q$	& \multicolumn{1}{|c} {$M= 28$} & \multicolumn{1}{|c|} {$M = 20$}
& \multicolumn{1}{|c|} {$M = 19$} & \multicolumn{1}{|c|} {$M = 18$} \\ \cline{2-5}
     	& \multicolumn{1}{|c} {$\v{E} = (4.1, 0.45)$} & \multicolumn{1}{|c|} {$\v{E} = (267, 21)$}
& \multicolumn{1}{|c|} {$\v{E} = (734, 56)$} & \multicolumn{1}{|c|} {$\v{E} = (2583, 464)$} \\ \hline
1/2       & $0.50\pm 0.08$   & $0.48\pm 0.08$ & $0.42\pm 0.07$ & $0.47\pm 0.08$\\
1/1       & $4.58\pm 0.23$   & $4.60\pm 0.23$ & $4.76\pm 0.23$ & $4.79\pm 0.23$\\
5/4       & $7.31\pm 0.29$   & $7.53\pm 0.29$ & $7.39\pm 0.29$ & $7.61\pm 0.29$\\
$\omega'$ & $71.65\pm 0.49$  & $71.06\pm 0.50$ & $71.02\pm 0.50$ & $70.73\pm 0.50$\\
3/2       & $12.05\pm 0.36$  & $12.35\pm 0.36$  & $12.18\pm 0.36$ & $12.07\pm 0.36$\\
2/1       & $2.72\pm 0.18$   & $2.72\pm 0.18$ & $2.87\pm 0.18$ & $3.05\pm 0.19$\\
5/2       & $0.97\pm 0.11$   & $0.99\pm 0.11$ & $1.08\pm 0.11$ & $1.03\pm 0.11$\\
3/1       & $0.22\pm 0.05$   & $0.27\pm 0.06$ & $0.28\pm 0.06$ & $0.24\pm 0.05$\\ \hline
\end{tabular}
\caption{Demonstration of the robustness of the computation of
attractor probabilities when $\g
= 10^{-5}$, $N = 18$ and the other parameters take the default values.
For comparison, the $M = 28$ column repeats the results in Table~\ref{vertab}.
The error estimates are $\v{E} = \max_{\v{x}_0\in\mathcal Q} (e_x, e_y)\times 10^{-14}$,
and are taken from Table~\ref{table1}. For all solutions and all values of
$M$, the 95\% confidence intervals overlap.}
\label{fuzztab}
\end{table}

\section{Discussion and Conclusions}
\label{conc_sec}

Weierstrass' Approximation Theorem~\cite{Hands} (real, multivariate 
polynomial version) states:

\begin{quote}
If $f$ is a continuous real-valued function defined on the set $[a,b]\times [c,d]$
and $\delta > 0$, then there exists a polynomial function $p$ in two variables such
that $|f(x,y) - p(x,y)| < \delta$ for all $x \in [a,b]$ and $y \in [c,d]$.
\end{quote}

In the light of this, it is not surprising that, for large enough degree $N$, and
number of timesteps per $2\pi$, $M$, the Frobenius method can give very good 
approximations to
the functions $X_i(\v{x})$, $Y_i(\v{x})$, $i = 1,\ldots, M$, that go to build up the
Poincar\'{e} map, $\v{P}(\v{x})$, and hence, to $\v{P}(\v{x})$ itself.
Less obvious is how effective a numerical ODE
solver based on such series approximations --- the high-order Euler method (HEM)
as we call it --- can be in practice.

In this paper, we have applied the HEM to a particular
problem, the spin-orbit problem, to illustrate its effectiveness in solving
this nonlinear ODE. We maintain that this is a non-trivial
problem, in the sense that the set of solutions can consist of many
coexisting periodic orbits as well as one quasi-periodic solution.
We show here that not only is the HEM capable of finding all the solutions
predicted by perturbation theory, and finds none that are not so predicted, but 
it also enables us to compute accurately the
mean frequency of the quasi-periodic solution and to make estimates of the probabilities
of the various coexisting attractors which agree with published results and the
Goldreich-Peale formula (where it applies). Additionally, compared to
standard numerical techniques, not only does HEM find all anticipated
solutions, but it is also about 40 times faster. All this is achieved by
using standard double precision arithmetic.

This increased speed comes at the cost of setting up the functions
$X_i(\v{x})$, $Y_i(\v{x})$, which map the solution and its derivative forward
by a time $h$, where $h = 2\pi/M$ is a timestep which is
relatively large since in practice $M\sim 25$. Before the advent of computer 
algebra this approach would have been impracticable for anything more than
very small $N$ --- in fact, too small, given our accuracy requirements ---
but such software is readily available nowadays and the process of setting 
up these functions is easily automated.

It would be useful to be able to define the class of ODEs for which HEM is
a good numerical algorithm. It not straightforward to define such a class,
although it is clear that all functions appearing in the ODE must
be sufficiently many times differentiable, in order for the Frobenius method
to work. When the dissipation is large, any standard ODE solver would
work; however, HEM comes into its own for problems with small dissipation.
Additionally, if the nonlinearity is multiplied by a small
parameter, conjecturally that may help to improve convergence --- see 
equations~(\ref{fourform}) and~(\ref{polyform}).

Further work is needed to investigate whether, for the spin-orbit problem
at least, the range of $y$-values can be extended. In fact, we have carried
out computations which show 
that the HEM works for at least $y\in[-5, 10]$, although at the cost of
increasing $M$ to 35--40. For $y$-values outside this range, we envisage 
that the polynomials $A_{i, j}(y)\ldots D_{i, j}(y)$ in equation~(\ref{polyform})
might be better replaced by a rational form, by using, for example, a 
Pad\'{e} approximation~\cite{Numrec}. 

Of interest too is the possibility that the technique described in~\cite{saari},
which uses a conformal transformation of the independent variable to extend
the radius of convergence of a power
series solution to an ODE, might be applied to the spin-orbit problem, thereby
allowing us to decrease $M$ and so further speed up our algorithm.
Another question for investigation is 
whether there exists a form in which the Poincar\'{e} map can be represented
that can be computed in significantly fewer arithmetical operations than
we have used here. 

Recently, it has been argued that MacDonald's model does not provide a
realistic description of the tidal torque and leads to
inconsistencies~\cite{makarov,will,noyelles}. In this paper, we have used
MacDonald's torque model in equation~(\ref{ode}) both for purposes of
comparison with existing literature --- in particular,~\cite{CC1,GP} --- and 
because its simplicity makes it particularly suitable for analytical
calculations. It would be interesting to investigate to what extent our
method could be applied to more general situations, such as those envisaged
in the papers quoted above.

\bigskip
{\large\bf Acknowledgements}
\medskip

The Authors wish to acknowledge the helpful comments received from the
anonymous reviewers, which have strengthened the paper by bringing to
our attention the freely
available software (DOP853 and the adaptive Taylor series method) with
which our algorithm has been compared in Section~\ref{speedsec}.
\null
\vfill

\appendix

\pagebreak

\centerline{\Large\bf Appendix: perturbation theory computations}

\section{Perturbation theory: periodic attractors}
\label{sec:2}

We carry out here the perturbation theory computations necessary to
establish thresholds for the periodic and quasi-periodic solutions which we
observe numerically. Details concerning this type of computation can be
found in~\cite{cod_lev,verhulst,gg_enc_a,gg_enc_b}.

\subsection{First order computation} \label{subsec:2.1}

We consider \eqref{ode} with $\gamma=\e \, C_1+\e^2C_2+O(\e^3)$ and
look for a solution $\v{x}(t)=(x(t),\,y(t))$ in the form of a
power series in $\e$, that is 
$\v{x}(t)=\v{x}^{(0)}(t)+\e\,\v{x}^{(1)}(t) + \e^{2}\v{x}^{(2)}(t)+\ldots$,
where $\v{x}^{(0)}(t)=( x_{0}+\omega_0 t, \omega_0)$, with $\omega_0=p/q$, and
$\v{x}^{(k)}(t)=(x^{(k)}(t),\,y^{(k)}(t))$ to be determined by requiring
that $\v{x}(t)$ be periodic in $t$ with period $2\pi q$.

A first order analysis gives
\begin{equation} \label{eq:2.2}
\begin{cases}
\dot  x^{(1)} = y^{(1)} , & \cr
\dot y^{(1)} = - G( x_{0}+\omega_0 t , t) - C_1 \alpha \left( \omega_0 - \omega \right) .
\end{cases}
\end{equation}
By introducing the Wronskian matrix
\begin{equation} \nonumber
W(t) = \left( \begin{matrix} 1 & t \\ 0 & 1  \end{matrix} \right) ,
\end{equation}
we can write $\v{x}^{(1)}(t)$ as
\begin{equation} \nonumber
\left( \begin{matrix}  x^{(1)}(t) \\ y^{(1)}(t) \end{matrix} \right) =
W(t) \left( \begin{matrix} \bar  x^{(1)} \\ \bar  y^{(1)} \end{matrix} \right) +
W(t) \int_{0}^{t} {\rm d}\tau \, W^{-1}(\tau)  \left( \begin{matrix} 0 \\ 
-G( x_{0}+\omega_0 \tau, \tau) - C_1 \alpha \left( \omega_0  - \omega \right)  \end{matrix} \right) ,
\end{equation}
with $\bar x^{(1)}=0$ and $\bar y^{(1)}$ to be fixed, so that
\begin{equation} \label{eq:2.3}
 x^{(1)}(t) = \bar  x^{(1)} + \bar y^{(1)} t -
\int_{0}^{t} {\rm d}\tau \int_{0}^{\tau} {\rm d}\tau' 
\left[ G( x_{0}+\omega_0 \tau', \tau') + C_1 \alpha \left( \omega_0 - \omega \right) \right] ,
\end{equation}
whereas $y^{(1)}(t)=\dot  x^{(1)}(t)$. For \eqref{eq:2.3} to be periodic
we have to require first of all that
\begin{equation} \label{eq:2.4}
M_{1}( x_{0}) : = \frac{1}{2\pi q} \int_{0}^{2\pi q} {\rm d} t
\left[ G( x_{0}+\omega_0 t, t) + C_1 \alpha \left( \omega_0 - \omega \right) \right] = 0 ,
\end{equation}
then fix $\bar y^{(1)}$ in such a way that the terms linear in $t$ in (\ref{eq:2.3}) cancel out.

Inserting \eqref{gdef} into \eqref{eq:2.4} leads to
\begin{equation} \nonumber
- \frac{1}{2\pi q} \sum_{k\in\calK} A_{k} \int_{0}^{2\pi q} {\rm d}t \,
\sin(2 x_{0}+2\omega_{0} t-kt) = C_1 \alpha \left( \omega_0 -\omega \right) 
\end{equation}
and hence
\begin{equation} \label{eq:2.5}
A_{k(p/q)} \sin 2 x_{0} = C_1 \alpha \Big( \frac{p}{q} - \omega \Big) ,
\qquad k(p/q) = \frac{2p}{q} .
\end{equation}
Since $A_{k} \neq 0$ only for $k\in\calK$, we have two possibilities:
\begin{enumerate}
\item if $\omega_0$ is of the form $\omega_0=p/2$, with $p\in\calK$, 
then for any $|C_{1}|<K_{1}(p)$, with
\begin{equation} \label{eq:2.6}
K_{1}(p) : = \frac{2|A_{p}|}{\alpha \, |p - 2\omega |} ,
\end{equation}
one can fix $ x_0$ in such a way that (\ref{eq:2.5}) is satisfied;
\item for all other values of $\omega_0$ one must require $C_1=0$.
\end{enumerate}

\subsection{Second order computation} \label{subsec:2.2}

The equations of motion to second order are
\begin{equation} \label{eq:2.7}
\begin{cases}
\dot  x^{(2)} = y^{(2)} , & \cr
\dot y^{(2)} = - \partial_{ x}G( x_{0}+\omega_0 t , t)\, x^{(1)}(t) - C_1 \alpha y^{(1)}(t) - C_2 \alpha \left( \omega_0 - \omega \right) ,
\end{cases}
\end{equation}
so that
\begin{equation} \label{eq:2.8}
 x^{(2)}(t) = \bar  x^{(2)} + \bar y^{(2)} t -
\int_{0}^{t} {\rm d}\tau \int_{0}^{\tau} {\rm d}\tau' 
\left[ \partial_{ x} G( x_{0}+\omega_0 \tau', \tau')\, x^{(1)}(\tau') + C_1 \alpha \dot x^{(1)}(\tau') +
C_2 \alpha \left( \omega_0 - \omega \right) \right] ,
\end{equation}
where
\begin{equation} \label{eq:2.9}
\partial_{ x} G( x, t) = \sum_{k\in\calK} 2 A_{k} \cos(2 x-kt) 
\end{equation}
and $ x^{(1)}(t)$ is obtained from (\ref{eq:2.3}). An explicit calculation gives
\begin{equation} \label{eq:2.10}
 x^{(1)}(t) = \bar x^{(1)} + 
\sum_{k\in\calK} \frac{A_k }{(2\omega_0-k)^2} \sin \left(2 x_0 + (2\omega_0-k)t \right) -
\sin (2 x_0) \sum_{k\in\calK} \frac{A_k }{(2\omega_0-k)^2} ,
\end{equation}
provided $\bar y^{(1)}$ is fixed in such a way that
\begin{equation} \nonumber
\bar y^{(1)} - \cos (2 x_0) \sum_{k\in\calK} \frac{A_k }{2\omega_0-k} = 0 . 
\end{equation}
Again for \eqref{eq:2.8} to be periodic we need that
\begin{equation} \label{eq:2.11}
M_2( x_{0}) : = \frac{1}{2\pi q} \int_{0}^{2\pi q} {\rm d} t
\left[ \partial_{ x} G( x_{0}+\omega_0 \tau', \tau')\, x^{(1)}(\tau') + C_1 \alpha \dot x^{(1)}(\tau') +
C_2 \alpha \left( \omega_0 - \omega \right) \right] = 0 ,
\end{equation}

If $\omega_0=p/2$, with $p\in\calK$, this simply produces a second order correction $\e^2C_{2}$ to the
leading order computed in the previous section. On the contrary, if $\omega_0$ is not of such form
then $C_1=0$ (by the analysis in Sect.~\ref{subsec:2.1}) and (\ref{eq:2.11}) becomes
\begin{equation} \label{eq:2.12}
M_2( x_{0}) : = \frac{1}{2\pi q} \int_{0}^{2\pi q} {\rm d} t
\left[ \partial_{ x} G( x_{0}+\omega_0 \tau', \tau')\, x^{(1)}(\tau') +
C_2 \alpha \left( \omega_0 - \omega \right) \right] = 0 ,
\end{equation}
By inserting (\ref{eq:2.10}) into (\ref{eq:2.12}) we find
\begin{equation} \nonumber
\frac{1}{2\pi q} \sum_{k,k'\in \calK} \frac{2 A_{k} A_{k'}}{(2\omega_0-k)^2}
\int_{0}^{2\pi q} {\rm d}t \,
\cos \left( 2 x_0 + (2\omega_0-k')t \right) \sin \left( 2 x_0 + (2\omega_0-k)t \right)=
C_{2} \alpha \left( \omega_0 - \omega \right) ,
\end{equation}
which implies
\begin{equation} \label{eq:2.13}
\sum_{\substack{k,k'\in \calK \\ k+k'=4\omega_0}} \frac{A_{k} A_{k'}}{(2\omega_0-k)^2}
\sin \left( 4  x_0 \right) =
C_{2} \alpha \left( \omega_0 - \omega \right) .
\end{equation}
Therefore, if $\omega_0$ is not of the form $\omega_0=p/2$, $p\in\calK$, two possibilities arise:
\begin{enumerate}
\item if $\omega_0$ is of the form $\omega_0=p/4$, with $p$ odd such that $p=k+k'\in \calK$, then,
defining
\begin{equation} \label{eq:2.14}
K_{2}(p) := \frac{16}{\alpha \, |p - 4 \omega|} \Biggl| \sum_{k=k_{1}(p)}^{k_{2}(p)}
\frac{A_{k} A_{p-k}}{(p-2k)^2} \Biggr|
\end{equation}
with $k_1(p)=\max\{-3,p-7\}$ and $k_2(p)=\min\{7,p+3\}$,
one has that for any $|C_2|<K_{2}(p)$ one can fix $ x_0$ in such a way that
(\ref{eq:2.13}) is satisfied;
\item for all other values of $\omega$ one must require $C_2=0$.
\end{enumerate}

\subsection{Threshold values} \label{subsec:2.4}

In the case of Mercury, one has $e=0.2056$ and hence $\bar L(e)=1.36937$,
$\bar N(e)=1.71971$, giving $\omega=1.25584$. For comparison, in the case
of the Moon, whose orbit is less eccentric, one has $e=0.0549$ and hence 
$\bar L(e)=1.02285$, $\bar N(e)=1.04135$, giving $\omega=1.01809$. 

Consider now the Sun-Mercury (S-M) system and, for comparison, the Earth-Moon (E-M) system.
To first order one finds the threshold values in Table~\ref{tab:5.1}, while to
second order the threshold values are as in Table~\ref{tab:5.2}.
For negative values of $p$, the values of the constants are less than $10^{-6}$ 
for E-M and less than $10^{-5}$ for S-M in Table \ref{tab:5.1},
less than $10^{-9}$ for E-M and less than $10^{-6}$ for S-M in Table \ref{tab:5.2}.

\medskip

\begin{table}[ht]
\centering
\begin{center}
\setlength\tabcolsep{5pt}
\vrule
\begin{tabular}{|l|l|l|l|l|l|l|l|} \hline
$p$ & 1  & 2  & 3 & 4  & 5 &  6  & 7 \\ \hline
E-M &  $\st5.178 \times 10^{-2}$ &
$\st53.64$ & $\st3.872\times 10^{-1}$ &
$\st2.533 \times 10^{-2}$ & 
$\st1.908 \times 10^{-3}$ &
$\st1.493 \times 10^{-4}$ &
$\st1.168 \times 10^{-5}$ \\ \hline
S-M & $\st 9.880 \times 10^{-2} $ &
$\st2.557$ &
$\st1.956$ &
$\st3.190 \times 10^{-1}$ &
$\st8.067 \times 10^{-2}$ &
$\st2.492 \times 10^{-2}$ &
$\st7.109 \times 10^{-3}$ \\ \hline
\end{tabular}
\hspace{-0.07cm}\vrule
\caption{Values of the constants $K_{1}(p)$ for $\alpha=\bar L(e)$ and $\omega=\nu(e)$.}
\label{tab:5.1}
\end{center}
\end{table}


\begin{table}[ht]
\centering
\begin{center}
\setlength\tabcolsep{4pt}
\vskip-.5truecm
\vrule
\begin{tabular}{|l|l|l|l|l|l|l|l|} \hline
$p$ & 1 & 3 & 5 & 7 & 9 & 11 & 13 \\ \hline
E-M & $\st3.909\times 10^{-6}$ &
$\st7.945\times 10^{-1}$ &
$\st6.386$ &
$\st5.531\times 10^{-2}$ &
$\st5.154 \times 10^{-4}$ &
$\st4.331 \times 10^{-6}$ &
$\st3.145 \times 10^{-8}$ \\ \hline
S-M & $\st1.200\times 10^{-4}$ &
$\st1.058$ &
$\st585.2$ &
$\st2.673$ &
$\st2.925\times 10^{-1}$ &
$\st3.507\times 10^{-2}$ &
$\st3.810\times 10^{-3}$ \\ \hline
\end{tabular}
\hspace{-0.08cm}\vrule
\caption{Values of the constants $K_{2}(p)$ for $\alpha=\bar L(e)$ and $\omega=\nu(e)$.}
\label{tab:5.2}
\end{center}
\end{table}

For $e=0.2056$ (Sun-Mercury system) and $\e=10^{-3}$, the threshold values corresponding 
to the resonances appearing in Tables~\ref{tab:5.1} and~\ref{tab:5.2} are given in 
Table~\ref{tab:5.3}.

\bigskip

\begin{table}[ht]
\centering
\begin{center}
\setlength\tabcolsep{3.5pt}
\vskip-0.5truecm
\hskip-0.4truecm
\vrule
\begin{tabular}{|l|l|l|l|l|l|l|l|} \hline
$\omega_0$ & 1/2 &  1 & 3/2 & 2 & 5/2 & 3 & 7/2 \\ \hline
\null & $\st9.880\times 10^{-5}$ &
$\st2.557 \times 10^{-3}$ &
$\st1.1956 \times 10^{-3}$ &
$\st3.190 \times 10^{-4}$ &
$\st8.067 \times 10^{-5}$ &
$\st2.492 \times 10^{-5} $ &
$\st7.109 \times 10^{-6}$ \\
\hline\noalign{\smallskip}\hline
$\omega_0$ & 1/4 & 3/4 & 5/4 & 7/4 & 9/4 & 11/4 & 13/4 \\ \hline
\null  & $\st1.200 \times 10^{-10}$ &
$\st1.058 \times 10^{-6}$ &
$\st5.852 \times 10^{-4}$ &
$\st2.673 \times 10^{-6}$ &
$\st2.926\times 10^{-7}$ &
$\st3.507\times 10^{-8}$ &
$\st3.810\times 10^{-9}$ \\ \hline
\end{tabular}
\hspace{-0.08cm}\vrule
\caption{S-M threshold values corresponding to the resonances listed in 
Tables~\ref{tab:5.1} and~\ref{tab:5.2} for $\e=10^{-3}$.}
\label{tab:5.3}
\end{center}
\end{table}

Therefore, for the system Sun-Mercury with $\e=10^{-3}$,
if $\gamma=10^{-5}$ the existing resonances are:
1:2, 1:1, 5:4, 3:2, 2:1, 5:2 and 3:1; if $\gamma=10^{-6}$ the existing resonances are
the same plus the further resonances 7:2, 3:4 and 7:4. 

\section{Perturbation theory: quasi-periodic attractors}
\label{sec:3}

We also look for a quasi-periodic solution of the form
\begin{equation} \label{eq:3.1}
 x(t) =  x_{0} + \omega' t + h( x_{0}+\omega' t,\e) , \qquad
h(\psi,\e) = \e \, h^{(1)}(\psi) + \e^{2} h^{(2)}(\psi) + \ldots 
\end{equation}
where $\omega'$ close to $\omega$ is to be determined.

The idea is to fix $\omega'$ and look for a solution of the form (\ref{eq:3.1}) to (\ref{ode})
with $\omega= \omega'+\mu(\omega',\e)$ for a suitable $\mu(\omega',\e)=\e \mu^{(1)}(\omega')+
\e^{2} \mu^{(2)}(\omega')+\ldots$. However, in (\ref{ode}) $\omega$ is a fixed parameter.
So, one should find the function $\mu(\omega',\e)$ and then try to solve the implicit function
problem $\omega'+\mu(\omega',\e)=\omega$. Unfortunately, the function
$\omega'\mapsto\mu(\omega',\e)$ is not smooth: a careful analysis shows that the function
is defined only for $\omega'$ satisfying a Diophantine condition. 
Nevertheless we do not address this problem here; we confine ourselves to
a third order analysis, neglecting any convergence problems; see 
Sect.~\ref{subsec:3.4} for further comments.

\subsection{First order computation} \label{subsec:3.1}

As in Sect.~\ref{sec:2} we write the differential equation~(\ref{ode}),
with $\omega=\omega'+\mu$ and $\gamma =C\e$, as an integral equation
\begin{equation} \label{eq:3.2}
 x(t) = \bar  x + \bar y \, t - \e \int_{0}^{t} {\rm d}\tau \int_{0}^{\tau} {\rm d}\tau'
\left[ G( x(\tau'),\tau') + C \alpha \bigl( \dot x(\tau')-\omega' - \mu \bigr) \right] ,
\end{equation}
We look for a solution of the form (\ref{eq:3.1}) and
set $ x^{(k)}(t)=h^{(k)}( x_{0}+\omega' t)$ for $k\ge 1$. 

Then to first order we obtain
\begin{equation} \label{eq:3.3}
 x^{(1)}(t) = \bar  x^{(1)} + \bar y^{(1)} t - \int_{0}^{t} {\rm d}\tau \int_{0}^{\tau} {\rm d}\tau'
\, G( x_{0}+\omega' \tau',\tau') ,
\end{equation}
which, after integration, gives
\begin{equation} \label{eq:3.4}
 x^{(1)}(t) = \tilde  x^{(1)} + \sum_{k\in\calK}
\tilde A_{k} \sin \left( 2 x_0 + (2\omega'-k)t \right) ,
\end{equation}
where
\begin{equation} \label{eq:3.5}
\tilde  x^{(1)} := \bar  x^{(1)} - \sin(2 x_0) \sum_{k\in\calK} \tilde A_{k} ,
\qquad \tilde A_{k} := \frac{A_{k}}{(2\omega'-k)^{2}} ,
\end{equation}
provided that $\bar y^{(1)}$ is fixed so as to satisfy
\begin{equation} \nonumber
\bar y^{(1)} - \cos(2 x_0) \sum_{k\in\calK} \frac{A_{k}}{2\omega'-k} = 0.
\end{equation}
%

\subsection{Second order computation} \label{subsec:3.2}

To second order (\ref{eq:3.2}) becomes
\begin{equation} \label{eq:3.6}
 x^{(2)}(t) = \bar  x^{(2)} + \bar y^{(2)} t - \int_{0}^{t} {\rm d}\tau \int_{0}^{\tau} {\rm d}\tau'
\left( \partial_{ x} G( x_{0}+\omega' \tau',\tau') \,  x^{(1)}(\tau') +
C \alpha \dot x^{(1)}(\tau') - C \alpha \mu^{(1)}(\omega') \right) .
\end{equation}
By using (\ref{eq:2.9}) and (\ref{eq:3.4}) we can write in (\ref{eq:3.6})
\begin{eqnarray}
&& \partial_{ x} G( x_{0}+\omega' t,t) \,  x^{(1)}(t) =
\sum_{k\in\calK} 2 A_{k} \tilde  x^{(1)} \cos (2 x_0+(2\omega'-k) t) \nonumber \\
& & \qquad \qquad + 
\sum_{k,k'\in\calK} 2 A_{k'} \tilde A_{k}
\cos (2 x_0+(2\omega'-k') t) \, \sin (2 x_0+(2\omega'-k) t) . \nonumber
\end{eqnarray}
Then, writing
\begin{equation} \nonumber
\cos (2 x_0+(2\omega'-k') t) \, \sin (2 x_0+(2\omega'-k) t) =
\frac{1}{2} \Bigl( \sin (4 x_0+(4\omega'-k-k') t) + \sin((k'-k)t) \Bigr) ,
\end{equation}
we find
\begin{eqnarray}
& & \int_{0}^{\tau} {\rm d}\tau'
\partial_{ x} G( x_{0}+\omega' \tau',\tau') \,  x^{(1)}(\tau') =
\sum_{k\in\calK} 2A_{k} \tilde  x^{(1)} \, \frac{\sin(2 x_0+(2\omega'-k)\tau)-
\sin(2 x_0)}{2\omega'-k} \nonumber \\
& & \qquad -
\sum_{k,k'\in\calK} A_{k'} \tilde A_{k} \, \frac{\cos(4 x_0+(4\omega'-k-k')\tau)-
\cos(4 x_0)}{4\omega'-k-k'} -
\sum_{\substack{k,k'\in\calK \\ k \neq k'}} A_{k'} \tilde A_{k} \, \frac{\cos((k'-k)\tau)-1}{k'-k} \nonumber
\end{eqnarray}
and hence in~(\ref{eq:3.6})
\begin{eqnarray}
& & - \int_{0}^{t} {\rm d} \tau \int_{0}^{\tau} {\rm d}\tau'
\partial_{ x} G( x_{0}+\omega' \tau',\tau') \,  x^{(1)}(\tau') =
\sum_{k\in\calK} 2 \tilde A_{k} \tilde  x^{(1)} \left( \cos(2 x_0+(2\omega'-k)t)-
\cos(2 x_0) \right) \nonumber \\
& & \qquad +
\sum_{k,k'\in\calK} A_{k'} \tilde A_{k} \, \frac{\sin(4 x_0+(4\omega'-k-k')t)-
\sin(4 x_0)}{(4\omega'-k-k')^2} +
\sum_{\substack{k,k'\in\calK \\ k \neq k'}} A_{k'} \tilde A_{k} \, \frac{\sin((k'-k)t)}{(k'-k)^{2}} \nonumber \\
& & \qquad + \, t \, \Biggl( \sum_{k\in\calK} \frac{2A_{k}\tilde x^{(1)} \sin(2 x_0)}{2\omega'-k} -
\sum_{k,k'\in\calK} A_{k'} \tilde A_{k} \, \frac{\cos(4 x_0)}{4\omega'-k-k'} -
\sum_{\substack{k,k'\in\calK \\ k \neq k'}} A_{k'} \tilde A_{k} \, \frac{1}{k'-k} \Biggr) . \nonumber
\end{eqnarray}
Furthermore in~(\ref{eq:3.6})
\begin{eqnarray} 
& & - \int_{0}^{t} {\rm d} \tau \int_{0}^{\tau} {\rm d}\tau' \, C \alpha \dot  x^{(1)}(\tau') =
- C \alpha \int_{0}^{t} {\rm d} \tau \left(   x^{(1)}(\tau) -  x^{(1)}(0) \right) \nonumber \\
& & \qquad = - \, C \alpha \sum_{k\in\calK} \tilde A_{k}
\int_{0}^{t} {\rm d} \tau \, \left(\sin(2 x_{0}+ (2\omega'-k)\tau) - \sin(2x_0)\right)\nonumber \\
& & \qquad = C \alpha \sum_{k\in\calK} \tilde A_{k} \,
\frac{\cos(2 x_{0}+(2\omega'-k)t) -\cos(2 x_0)}{2\omega'-k}
+ t\, C\,\alpha \sin(2 x_0)\sum_{k\in\calK} \tilde A_{k}  , \nonumber
\end{eqnarray}
where (\ref{eq:3.4}) has been used. The coefficient $\mu^{(1)}(\omega')$ in (\ref{eq:3.6})
has to be fixed so as to cancel out any term linear in $\tau$ produced by the $\tau'$-integration,
if such a term exists. Since there is no such term, we set
$\mu^{(1)}(\omega')=0$. Therefore, if we also set
\begin{equation} \nonumber
\bar y^{(2)} + \sin(2 x_0) \sum_{k\in\calK} \frac{2A_{k}\tilde x^{(1)} }{2\omega'-k} -
\cos(4 x_0) \sum_{k,k'\in\calK} \frac{A_{k'} \tilde A_{k}}{4\omega'-k-k'} -
\sum_{\substack{k,k'\in\calK \\ k \neq k'}} \frac{A_{k'} \tilde A_{k}}{k'-k}
+ C \alpha \sin(2 x_0)\sum_{k\in\calK} \tilde A_{k} = 0 ,
\end{equation}
we obtain
\begin{eqnarray} \label{eq:3.7}
& &  x^{(2)}(t) = \tilde x^{(2)} + \sum_{k\in\calK} \tilde B_{k} \cos(2 x_0+(2\omega'-k)t) \nonumber \\
& & \qquad \qquad +
\sum_{k,k'\in\calK} \tilde C_{k,k'} \sin(4 x_0+(4\omega'-k-k')t) 
+ \sum_{\substack{k,k'\in\calK \\ k \neq k'}} \tilde D_{k,k'} \sin((k'-k)t) ,
\end{eqnarray}
where we have defined
\begin{eqnarray} \label{eq:3.8}
& & \tilde x^{(2)} = \bar  x^{(2)} -
\cos(2 x_0) \sum_{k\in\calK} \tilde B_{k} -
\sin(4 x_0) \sum_{k,k'\in\calK} \tilde C_{k,k'} , \nonumber \\
& & \tilde B_{k} := 
2 \tilde A_{k} \tilde  x^{(1)} + \frac{C \alpha \tilde A_{k}}{2\omega'-k} , \qquad
\tilde C_{k,k'} := \frac{A_{k'} \tilde A_{k}}{(4\omega'-k-k')^2} , \qquad
\tilde D_{k,k'}:= \frac{A_{k'} \tilde A_{k}}{(k'-k)^{2}} . 
\end{eqnarray}

\subsection{Third order computation} \label{subsec:3.3}

To third order we have
\begin{eqnarray} \label{eq:3.9}
& &  x^{(3)}(t) = \bar  x^{(3)} + \bar y^{(3)} t - \int_{0}^{t} {\rm d}\tau \int_{0}^{\tau} {\rm d}\tau'
\Bigl( \partial_{ x} G( x_{0}+\omega' \tau',\tau') \,  x^{(2)}(\tau') \nonumber \\
& & \qquad \qquad \qquad \, + \; \frac{1}{2} \partial_{ x}^{2} G( x_{0}+\omega' \tau',\tau') \, ( x^{(1)}(\tau'))^{2}
+ C \alpha \dot x^{(2)}(\tau') - C \alpha \mu^{(2)}(\omega') \Bigr) ,
\end{eqnarray}
where once more $\mu^{(2)}(\omega')$ has to be fixed in such a way that the $\tau'$-integration
does not produce any term linear in $\tau$.

If we only want to determine $\mu(\omega', \epsilon)$ to second order, then we do not 
need to compute $ x^{(3)}(t)$ --- which would be needed to compute $\mu^{(3)}(\omega')$
--- and we have only to single out the terms linear in $\tau$ arising from
\begin{equation} \label{eq:3.10}
\int_{0}^{\tau} {\rm d}\tau'
\Bigl( \partial_{ x} G( x_{0}+\omega' \tau',\tau') \,  x^{(2)}(\tau')
+ \frac{1}{2} \partial_{ x}^{2} G( x_{0}+\omega' \tau',\tau') \, ( x^{(1)}(\tau'))^{2} \Bigr) ,
\end{equation}
where we have also used the fact that no term linear in $\tau$ is produced by the 
integration of $C \dot x^{(2)}(\tau')$.

We have in~(\ref{eq:3.9})
\begin{eqnarray} \label{eq:3.11}
& & \partial_{ x} G( x_{0}+\omega' t,t) \,  x^{(2)}(t)
+ \frac{1}{2} \partial_{ x}^{2} G( x_{0}+\omega' t,t) \, ( x^{(1)}(t))^{2} =
\sum_{k\in\calK} 2A_{k}\tilde x^{(2)} \cos (2 x_{0}+(2\omega' -k)t) \nonumber \\
& & \qquad + \; 
\sum_{k,k'\in\calK} 2A_{k'}\tilde B_{k} \cos (2 x_{0}+(2\omega' -k')t) \cos (2 x_{0}+(2\omega' -k)t) \nonumber \\
& & \qquad + \; 
\sum_{k,k',k''\in\calK} 2A_{k''}\tilde C_{k,k'} \cos (2 x_{0}+(2\omega' -k'')t) \sin (4 x_{0}+(4\omega' -k-k')t) \nonumber \\
& & \qquad + \; 
\sum_{\substack{k,k',k''\in\calK \\ k\neq k'}}
2A_{k''}\tilde D_{k,k'} \cos (2 x_{0}+(2\omega' -k'')t) \sin ((k'-k)t)  \\
& & \qquad - \; 
\sum_{k\in\calK} 2A_{k} (\tilde x^{(1)})^{2} \sin (2 x_{0}+(2\omega' -k)t) \nonumber \\
& & \qquad - \; 
\sum_{k,k'\in\calK} 4A_{k'} \tilde  x^{(1)}\tilde A_{k} \sin (2 x_{0}+(2\omega' -k')t) \sin (2 x_{0}+(2\omega' -k)t) \nonumber \\
& & \qquad - \; 
\sum_{k,k',k''\in\calK} 2A_{k''}\tilde A_{k} \tilde A_{k'} \sin (2 x_{0}+(2\omega' -k'')t)
\sin (2 x_{0}+(2\omega' -k)t) \sin (2 x_{0}+(2\omega' -k')t) , \nonumber 
\end{eqnarray}
where both (\ref{eq:3.4}) and (\ref{eq:3.7}) have been used.

If we use the trigonometric identities
\begin{eqnarray}
& & \cos\alpha \cos \beta = \frac{1}{2} \Bigl( \cos(\alpha+\beta) + \cos(\alpha-\beta) \Bigr) , 
\;\;\;\; \cos\alpha \sin \beta = \frac{1}{2} \Bigl( \sin(\alpha+\beta) + \sin(\beta-\alpha) \Bigr) , \nonumber \\
& & \sin\alpha \sin \beta = \frac{1}{2} \Bigl( \cos(\alpha-\beta) - \cos(\alpha+\beta) \Bigr) , \nonumber \\
& & \sin\alpha \sin \beta \sin \gamma = \frac{1}{4} \Bigl( \sin(\alpha-\beta+\gamma) +
\sin(\gamma-\alpha+\beta) - \sin(\alpha+\beta+\gamma) - \sin(\gamma-\alpha-\beta)  \Bigr) , \nonumber 
\end{eqnarray}
we realise immediately that only the second and sixth lines in (\ref{eq:3.11}) produce terms linear in $\tau$ after integration.
Indeed one has in~(\ref{eq:3.11})
\begin{eqnarray} \label{eq:3.11bis}
& & \sum_{k,k'\in\calK} 2A_{k'} \tilde B_{k} \cos (2 x_{0}+(2\omega' -k')t) \cos (2 x_{0}+(2\omega' -k)t) \nonumber \\
& & \qquad \qquad = \; 
\sum_{k,k'\in\calK} A_{k'} \tilde B_{k}  \cos (4 x_{0}+(4\omega' -k-k')t) 
+ \sum_{k,k'\in\calK} A_{k'} \tilde B_{k} \cos ((k' -k)t) ,
\end{eqnarray}
so that the term with $k=k'$ in the second sum in~(\ref{eq:3.11bis}) gives
\begin{equation} \label{eq:3.12}
\int_{0}^{\tau} {\rm d}\tau' 
\sum_{k\in\calK} A_{k} \tilde B_{k} = \tau \sum_{k\in\calK} A_{k} \tilde B_{k} .
\end{equation}
and, analogously, in~(\ref{eq:3.11}), one has
\begin{eqnarray} \label{eq:3.12bis}
& & - \sum_{k,k'\in\calK} 4A_{k'} \tilde  x^{(1)}\tilde A_{k} \sin (2 x_{0}+(2\omega' -k')t)
\sin (2 x_{0}+(2\omega' -k)t) \nonumber \\
& & \qquad \qquad = \; 
- \sum_{k,k'\in\calK} 2A_{k'} \tilde  x^{(1)}\tilde A_{k} \cos ((k'-k)t) +
\sum_{k,k'\in\calK} 2A_{k'} \tilde  x^{(1)}\tilde A_{k} \cos (4 x_{0}+(4\omega' -k-k')t) ,
\end{eqnarray}
so that the term with $k=k'$ in the first sum in~(\ref{eq:3.12bis}) gives 
\begin{equation} \label{eq:3.13}
\int_{0}^{\tau} {\rm d}\tau' \Biggl( - 
\sum_{k\in\calK} 2A_{k} \tilde  x^{(1)}\tilde A_{k} \Biggr) =
- \tau \sum_{k\in\calK} 2 A_{k} \tilde  x^{(1)}\tilde A_{k} .
\end{equation}
By collecting together the contributions (\ref{eq:3.12}) and (\ref{eq:3.13})
with that arising from the term in $\mu^{(2)}(\omega')$ in (\ref{eq:3.9}), we find
\begin{equation} \nonumber
\tau \Biggl( \; \sum_{k\in\calK} \bigl( A_{k} \tilde B_{k} - 2 A_{k} \tilde  x^{(1)}\tilde A_{k} \bigr) -
C \alpha \mu^{(2)}(\omega') \Biggr) .
\end{equation}
By (\ref{eq:3.8}) we have
\begin{equation} \nonumber
A_{k} \tilde B_{k} - 2 A_{k} \tilde  x^{(1)}\tilde A_{k} =
2 A_{k} \tilde A_{k} \tilde  x^{(1)} + \frac{C \alpha A_{k}\tilde A_k}{2\omega'-k} -
2 A_{k} \tilde  x^{(1)}\tilde A_{k} = \frac{C \alpha A_{k}\tilde A_k}{2\omega'-k} .
\end{equation}
so that one has to fix
\begin{equation} \label{eq:3.14}
\mu^{(2)}(\omega') = \sum_{k\in \calK} \frac{A_{k}\tilde A_k}{2\omega'-k}
= \sum_{k\in \calK} \frac{A_k^2}{(2\omega'-k)^3} .
\end{equation}
An explicit computation gives, in the case of Mercury, $\mu^{(2)}(\omega)=2.284502$ and,
in the case of the Moon, $\mu^{(2)}(\omega)=7.040139$.

\subsection{Conclusions} \label{subsec:3.4}

By requiring $\om'$ to satisfy a Diophantine condition such as
\begin{equation} \label{eq:3.15}
\left| \om'\nu_1+\nu_2 \right| \ge
\frac{\gamma_{0}}{(|\nu_1|+|\nu_2|)^{\tau_{0}}} ,
\end{equation}
where $\nu_1$ and $\nu_2$ are integers, and with $\gamma_0>0$ and $\tau_0>1$,
the analysis can be pushed to any perturbation order. The series for
$\mu(\om',\e)$
can then be proved to converge to a function
$\mu(\omega',\e)=\e^{2}\mu^{(2)}(\om')+O(\e^{3})$
depending analytically on $\e$.
In fact, this has been proved in \cite{CC2} --- it could also be proved
directly, by using diagrammatic
techniques (see for instance~\cite{G} for a review) to show
that, to any perturbation order $k$, the functions $x^{(k)}(t)$ and the
coefficients $\mu^{(k)}(\om')$
are bounded above proportionally to a constant to the power $k$.
Moreover, both the solution (\ref{eq:3.1}) and the
function $\mu(\om',\e)$ are not smooth in $\omega'$: in fact they are
defined on a Cantor set $\Omega$.
However, the function admits a Whitney extension \cite{W,CG,P} to a $C^{\io}$
function, so that one can consider the implicit function problem
\begin{equation} \label{eq:3.16}
\omega'+\mu(\omega',\e)=\omega .
\end{equation}
Such an equation admits a solution
\begin{equation}
\omega' = \omega - \e^{2} \mu^{(2)}(\omega') + O(\e^{3}) ,
\label{omp}
\end{equation}
so that, if for a fixed $\e$ the corresponding $\omega'$ is Diophantine,
then we
have a quasi-periodic attractor of the form (\ref{eq:3.1}).

For $\e_0>0$ and $\omega$ Diophantine, the set of values
$\e\in[0,\e_{0}]$
such that $\omega'$ is Diophantine has full measure in $[0,\e_{0}]$.
However the convergence of the series
requires for $\e\gamma_0^2$ to be small, so that the set of values of $\e$
for which the quasi-periodic
attractor exists has large, but not full measure. So, for fixed $\omega$,
it is a non-trivial problem
to understand whether a smooth quasi-periodic attractor can exist. Indeed, for
fixed $\omega$ and $\e$ one has
first to compute the solution $\om'$ to the implicit equation
(\ref{eq:3.16}) and then to check whether
such a solution satisfies the Diophantine condition (\ref{eq:3.15}).

\section{Computation of \mbox{\boldmath $\Delta\omega$} for 
\mbox{\boldmath $\varepsilon = 10^{-6}$}}

For $\varepsilon = 10^{-6}$, double precision arithmetic is
inadequate to estimate $\Delta\omega$: in this case, we are after all
attempting to find a difference of order $10^{-12}$ between two numbers,
$\omega$ and $\omega'$, both of order unity. Furthermore, this difference can only 
be estimated by iterating many times a (HEM-approximated) Poincar\'{e} map,
in which the error per iteration is $O(10^{-14})$.  In fact, we estimate
$\Delta\omega\approx 5.7\times 10^{-12}$ using double precision arithmetic
and $10^8$ iterations. Thus, use of a higher-accuracy computation is indicated.

We therefore use a HEM  with $M = 25, N = 20$, for which the maximum value
of $e_x \approx 2.6\times 10^{-21}$ when we iterate it using CA with 35
significant figures. Since this accurate CA implementation is about 8000
times slower that the equivalent LL computation, we reduce the number of
iterations $n$ to $10^6$ in equation~(\ref{omp_comp}). Also, convergence to
$\omega'$ is quite slow, so we extrapolate to estimate the limit as
$n\rightarrow\infty$.

\begin{figure}[!ht]
\centering
\includegraphics*[width=3in]{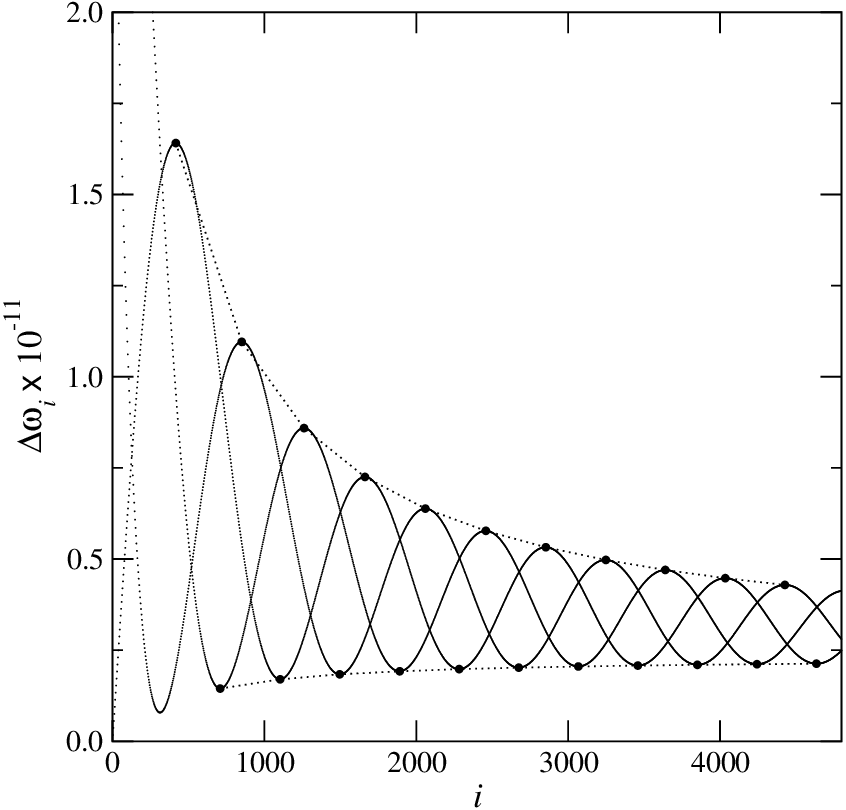}
\caption{A plot of $\Delta\omega_i$, defined in the text, against $i$,
showing convergence to the asymptotic value, $\Delta\omega$. The dashed
lines show the least squares fit curves through the peaks and troughs of the plot of
$\Delta\omega_i$.}
\label{interp}
\end{figure}

In order to illustrate this convergence and extrapolation, we include
Fig.~\ref{interp}, which is a plot of
$\Delta\omega_i = [x(2\pi i B) - x(0)]/(2\pi i B)$ against $i$ with $B = 200$.
Superimposed on the plot are the peak and trough values, shown as filled
circles, and a least squares fit curve (dashed line) through just these values. The curve
is of the form $y = a_0 + a_1/i + a_2/i^2 + a_3/i^3$ and for the peak
values, $a_0 = 2.278\times 10^{-12}$; for the trough values, $a_0 =
2.284\times 10^{-12}$. We therefore estimate $\Delta\omega \approx
2.28\times 10^{-12}$ for $\varepsilon = 10^{-6}$. This value should be
compared with that given by perturbation theory, which is $2.284\times
10^{-12}$.

\section{The Fourier form of the expressions for
\mbox{\boldmath $X_i(\v{x}), \, Y_i(\v{x})$} }

We start by expanding the sine terms in $G(x, t)$ in equation~(\ref{gdef}) in the 
spin-orbit ODE to obtain
\begin{equation} \label{ode_spec}
\begin{cases}
\dot  x = y , \\
\dot y = -\e[a(t)\cos 2x + b(t)\sin 2x] -\gamma\alpha\left(y - \omega\right)
\end{cases}
\end{equation}
where $a(t),\, b(t)$ are Fourier polynomials in $t$. We wish to
show that, formally, we can write
\begin{equation}
x(t_0 + h) = x_0 + A_{i, 0}(y_0) + 
\sum_{j = 1}^\infty \e^j\left[A_{i, j}(y_0)\cos 2jx_0 + 
B_{i, j}(y_0)\sin 2jx_0\right]
\label{eff}
\end{equation}
where $x_0 = x(t_0)$, $y_0 = y(t_0)$ and $A_{i, j},\, B_{i, j}$ are
polynomials in $y_0$, $h$ and the parameters in the ODE, but not $x_0$.
We refer to this as the Fourier series form.  The point we wish to make here is 
that the $j$-th coefficient in the expansion of $x(t_0 + h)$ in this form always 
has a factor of $\e^j$. That is not to say that, for small $\e$, the Fourier 
coefficients themselves decrease exponentially with
increasing $j$, because we do not prove that $A_{i, j},\, B_{i, j}$ grow
more slowly than exponentially.

To show that $x(t_0 + h)$ can be written in the form~(\ref{eff}), we start
with the Taylor series expansion
$$x(t_0 + h) = x_0 + h y_0 + \sum_{i = 2}^\infty \frac{h^i}{i!} y_0^{(i-1)},$$
where $y_0^{(i)}$ is the $i$-th derivative of $y(t) = \dot{x}(t)$ at $t=t_0$.
Using the fact
that the ODE~(\ref{ode_spec}) supplies us with a means for substituting for the
first --- and hence, recursively, for all --- derivatives of $y$, then
$y^{(i)}(t)$ can be written as a sum of terms each of which is a 
product of the form
$$T(t) = \nu' p(t) y(t)^k c^{d-l} s^l,$$
where $\nu'$ is a numerical constant; $p(t)$ is a combination of $a(t)$
and $b(t)$ and their derivatives; integer $k\geq 0$; 
$c = \cos(2x(t)),\, s = \sin(2x(t))$;
and $d\geq l \geq 0$. We wish to prove that all terms have a factor
of $\e^d$, so that we can always write $\nu' = \nu\e^d$, where $\nu$
is another constant; from this, equation~(\ref{eff}) will follow.

Let us define the degree of the term $T$ as the integer $d$, so that the degree
of $T$ is the sum of the powers of $c$ and $s$ appearing in $T$.

Differentiating $T$ with respect to $t$, we obtain
$$(\nu')^{-1}\,\dot{T}(t) = \dot{p} y^kc^{d-l}s^l + 
2 p y^{k+1}\left[-(d-l) c^{d-l-1} s^{l+1} + l c^{d-l+1} s^{l-1}\right]
+  k py^{k-1}\dot{y} c^{d-l} s^l.$$
As it stands, this expression consists of four terms each of degree $d$,
but using~(\ref{ode_spec}) to replace $\dot{y}$, the last term becomes
$$kp y^{k-1}\left[-\e a c - \e b s - \gamma\alpha y + \gamma\alpha\omega\right]c^{d-l} s^l = 
-\e kpy^{k-1}\left[a c^{d-l+1}s^l + b c^{d-l} s^{l+1}\right]
+kp y^{k-1}\gamma\alpha(-y+\omega)c^{d-l} s^l. $$
This expression consists of two terms of degree $d+1$, both of which are
multiplied by $\e$, and two terms of degree $d$, neither of which are 
multiplied by $\e$. Hence, differentiation of a term of degree $d$, followed by
substitution of $\dot{y}$, if present, leads to an expression of the form
$\e\times[\mbox{sum of terms of degree $(d+1)$}] + [\mbox{sum of terms of degree $d$}]$.
Since differentiation and substitution are the only processes by which
$y^{(i)}$ is generated, by induction all terms in $y^{(i)}$ of degree $d$ have a
factor of $\e^d$. Furthermore, any expression $c^{d-l} s^l$ is equal to a
sum of terms of the form $\sin 2ix$, $\cos 2ix$, with $i = 0, \ldots, d$;
and from this, the form~(\ref{eff}) follows.

\begin{figure}[!ht]
\centering
\includegraphics*[width=3in]{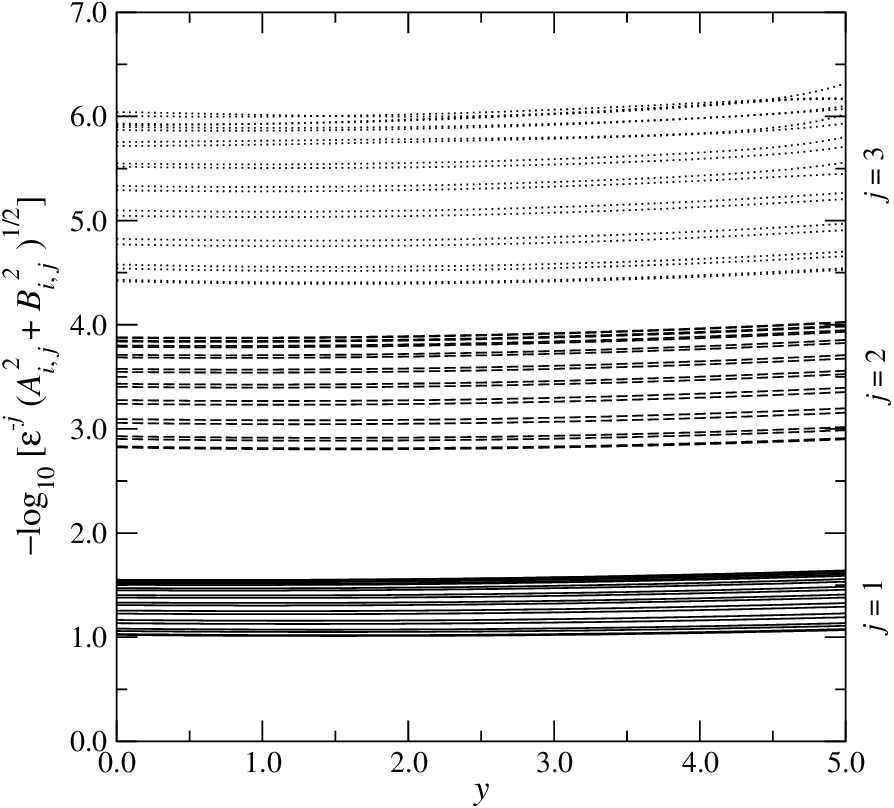}
\caption{A logarithmic plot of $\e^{-j}\sqrt{A_{i, j}(y)^2 + B_{i,
j}(y)^2}$, $i = 1,\ldots M = 20$ and $ j = 1, 2, 3$, against $y$,
where the polynomials $A_{i, j}$ and $B_{i, j}$ are defined in equation~(\ref{eff}).
Only the first three Fourier coefficients are needed to meet the error
criterion explained in Sect.~\ref{setup_sec}.
The polynomials were computed for $\epsilon = 1.2\times 10^{-4}$, $K =
10^{-4}$ and $e = 0.2056$. The figure shows that, for these parameters at
least, $A_{i, j}$ and $B_{i, j}$ decrease faster than than $\e^j$ for
all $i$. }
\label{fform}
\end{figure}

Numerical evidence suggests that $A_{i, j}$ and $B_{i, j}$ actually
decrease faster than $\e^j$, at least for $j = 2,\, 3$ --- see Fig.~\ref{fform}.

\vfill
\pagebreak

\end{document}